# BAYESIAN ANALYSIS IN MOMENT INEQUALITY MODELS

By Yuan Liao and Wenxin Jiang

*Northwestern University*

This paper presents a study of the large-sample behavior of the posterior distribution of a structural parameter which is partially identified by moment inequalities. The posterior density is derived based on the limited information likelihood. The posterior distribution converges to zero exponentially fast on any $\delta$-contraction outside the identified region. Inside, it is bounded below by a positive constant if the identified region is assumed to have a nonempty interior. Our simulation evidence indicates that the Bayesian approach has advantages over frequentist methods, in the sense that, with a proper choice of the prior, the posterior provides more information about the true parameter inside the identified region. We also address the problem of moment and model selection. Our optimality criterion is the maximum posterior procedure and we show that, asymptotically, it selects the true moment/model combination with the most moment inequalities and the simplest model.

## 1. Introduction.

1.1. *Formulation of the problems.* Let $(\bar{\Omega}, \mathcal{A}, P)$ denote a probability space. Suppose that we are interested in some structural parameter $\theta_0 \in \mathbb{R}^d$ that satisfies a set of *moment inequality* conditions:

$$Em_j(X, \theta_0) \geq 0, \qquad j = 1, \ldots, p, \tag{1.1}$$

where $m_j(\cdot, \theta), i = 1, \ldots, p$, are known real-valued moment functions. $X$ is an observable random vector defined on $(\bar{\Omega}, \mathcal{A}, P)$ and we assume that we observe i.i.d. or stationary realizations $X^n = \{X_1, \ldots, X_n\}$ of $X$. A model that is characterized by moment inequalities (1.1) is usually called a *moment inequality model*.

A key feature of moment inequality models is that $\theta_0$ is not necessarily point identified: there exists more than one solution to the inequalities in (1.1) if $Em_j(X, \theta_0)$ is viewed as a function of $\theta_0$. In other words, if we let $\Theta$ be the parameter space that contains $\theta_0$ and define

(1.2) $$\Omega = \{\theta \in \Theta, Em_j(X, \theta) \geq 0, j = 1, \ldots, p\},$$

then $\Omega$ can be a nonsingleton set. In this case, we say that $\theta_0$ is *partially identified* on $\Omega$ and $\Omega$ is called the *identified region*.

Many partially identified models are characterized by such moment inequalities, where the parameter of interest is only partially identified and therefore cannot possibly be consistently estimated. In this framework, since the identified region captures all of the information about the parameter, it becomes one of the most interesting subjects of study in moment inequality models [see, e.g., Chernozhukov, Hong and Tamer (2007), CHT, hereafter].

In addition to the problem of studying the identified region, there is also a *moment/model selection* problem in moment inequality models. Suppose that we have $p$ candidate moment inequalities

$$Em_j(X, \theta) \geq 0, \qquad j = 1, \ldots, p,$$

with a $k$-dimensional parameter vector $\theta = (\theta_1, \ldots, \theta_k)^T$ that belongs to the parameter space $\Theta_1 \times \cdots \times \Theta_k$. The moment selection problem refers to selecting the best subset of the moment inequalities among all of the possible candidates, while the model selection procedure addresses the problem of selecting the best model that is characterized by setting some components of the parameter to be zero. Such a candidate model can be a parameter subspace like $\{0\} \times \Theta_2 \times \cdots \times \Theta_k$. Therefore, the moment/model selection procedure produces a combination of moment inequalities and a parameter subspace. For instance, in Example 1.3 regarding the instrumental variable regression with interval censoring model, the moment selection problem can correspond to selecting the instrumental variables (components of $Z$), while the model selection problem is related to selecting the useful explanatory variables (components of $X$) that have nonzero regression coefficients. Ultimately, the selected combination should achieve some sense of optimality.

1.2. *Some motivating examples.* There are several interesting examples for the moment inequality models described above, where the parameter of interest is identified on a nonsingleton set.

EXAMPLE 1.1 [Interval censored data; see, e.g., Example 1 of CHT (2007)]. Let $Y$ be a real-valued random variable which lies in $[Y_1, Y_2]$ almost surely; $Y_1$ and $Y_2$ are observed random variables, but $Y$ is not observed. (Sometimes one may assume that $Y_2 = Y_1 + 1$, as in the case where $Y_1$ is the recorded integer part of $Y$.) The parameter of interest is $\theta_0 = E(Y)$. We then have the following moment inequalities:

$$E(Y_2 - \theta_0) \geq 0, \qquad E(\theta_0 - Y_1) \geq 0.$$

Then $\theta_0$ is partially identified on $\Omega = [EY_1, EY_2]$.

EXAMPLE 1.2 [Missing data; see, e.g., Example 1 of Canay (2008)]. Assume that $(Y, Z) \in [0,1] \times \{0,1\}$ and that $Y$ is observed only when $Z = 1$. Suppose that we are interested in the parameter $\theta_0 = EY$, which corresponds to $E[ZY + (1-Z)Y - \theta_0] = 0$, where $ZY$ is known, but not $(1-Z)Y$. Noting that $0 \leq (1-Z)Y \leq 1-Z$, we have moment inequalities

$$E(\theta_0 - ZY) \geq 0, \qquad E(ZY - \theta_0 + 1 - Z) \geq 0.$$

Then $\theta_0$ is partially identified on $\Omega = [E(ZY), E(ZY) + 1 - EZ]$.

EXAMPLE 1.3 [Interval regression model; see, e.g., Example 2 of CHT (2007)]. Consider $Y, Y_1, Y_2$, as in the setup of Example 1.1, assuming that the conditional mean of the unobserved $Y$ is $X^T \theta_0$, where $\theta_0$ is the parameter of interest and $X$ is a regressor vector. Due to $Y_1 \leq Y \leq Y_2$, we then have moment inequalities

$$EZ(Y_2 - X^T\theta_0) \geq 0, \qquad EZ(X^T\theta_0 - Y_1) \geq 0,$$

where $Z$ is a vector of positive functions of $X$ or positive instrumental variables.

1.3. *Literature review and contributions of this paper.* Many frequentist inference procedures for the identified region as well as the true parameter have been developed in this growing area of interest. For example, Chernozhukov, Hong and Tamer (2007) construct an econometric criterion function so that its set of minimizers form the identified region. They consider consistent estimation of the identified region and have shown that their set estimator is consistent in Hausdorff distance. Additionally, they derive the convergence rate of their estimator and construct the confidence set for the identified region. Moreover, Andrews and Soares (2007) develop confidence sets of the identified region and a test of the moment inequalities/equalities based on generalized moment selection. Among others, Rosen (2008) provides a formulation of criterion functions that differ from CHT and derives analytical critical values of the confidence region. Beresteanu and Molinari (2008) recently proposed inference procedures when the identified region can be written as a transformation of the Aumann expectation based on random set theory. Some additional papers in the literature that consider inference with partially identified models include Pakes et al. (2006), Andrews and Jia (2008), Romano and Shaikh (2008), Bugni (2007), Horowitz and Manski (2000), Manski and Tamer (2002), Canay (2008) and Liu and Shao (2003).

This paper studies a Bayesian approach to the moment inequality models. The Bayesian procedure provides distributional information for the partially identified parameter, both inside and outside the identified region, through its posterior distribution. The advantages of using posterior distributions to characterize the parameters are many. First, as pointed out by Poirier (1998), a Bayesian analysis of partial identification models is always possible if a proper prior for the parameters is specified. If we have some a priori information on $\theta_0$, then, by using a properly chosen prior distribution, the resulting posterior density may not be flat within the identified region; this provides evidence that the parameter is more likely to lie in some particular area. Second, even with a flat prior distribution, when $\theta_0$ is multidimensional, the posterior density of some components of $\theta_0$ may no longer be flat, due to the shape of the identified region. Hence, if we are interested in these components of $\theta_0$, then the posterior density can still provide extra information on their locations within the identified region. As a third advantage, it can be shown asymptotically that the posterior density has support only on the identified region. Containing more information, a posterior density can always be used to estimate the identified region, but not vice versa. Finally, the MCMC method is a very powerful method to draw samples from the posterior, which can be used for approximations of the calculation of the posterior statistics. In addition, those posterior samples can also be used in frequentist methods to estimate the identified region, by, for example, minimizing an econometric criterion function in CHT.

In fact, Bayesian methods have been extensively applied to nonidentified situations. Gelfand and Sahu (1999) have studied issues surrounding nonidentifiability and improper priors in the context of generalized linear models. Neath and Samaniego (1997) consider Bayesian updating for a nonidentified two-parameter binomial model. Gustafson (2005) studies Bayesian inference in nonidentified scenarios involving misclassification and measurement errors, which was discussed by a number of prominent researchers. Recently, Moon and Schorfheide [(2009), hereafter, MS] have considered the Bayesian approach to partially identified models when the model can involve three types of parameters: the structural parameters of interest, a reduced-form parameter vector that is point-identified by data and also a vector of auxiliary parameters which links the structural and reduced-form parameters via some known function. They also derive the Bayesian credible sets and compare them with frequentist confidence intervals for a number of particular models. All of these papers use traditional posteriors based on the likelihood function instead of the moment inequalities.

Our Bayesian approach proceeds within a more general framework. In contrast to the previous work, we do not need to have a full probability model for the observed data. Starting from moment inequalities $Em(X, \theta_0) \geq 0$, where $m(X, \cdot)$ is a known function of $\theta_0$, we introduce a bias parameter $\lambda_0 \geq 0$ so that $Em(X, \theta_0) = \lambda_0$ and place prior distributions on $(\theta_0, \lambda_0)$. The posterior density of $\theta_0$ can then be derived based on a limited information likelihood function, which is generated by the conditional asymptotic distribution of $\frac{1}{n}\sum_{i=1}^{n} m(X_i, \theta_0) - \lambda_0$ given $(\theta_0, \lambda_0)$, integrating out $\lambda_0$. We study in detail the frequentist behavior of the posterior density function of $\theta_0$. We derive the bounds of convergence rates of the posterior density, both inside and outside of the identified region. We show that there is a large "gap" between them. Once the posterior density and its frequentist properties are obtained, it is easy to derive consistent estimators for the identified region. However, we point out that a posterior density provides more information than a region estimation since it can also incorporate prior information and describe how likely the true parameter is to be distributed both inside and outside the identified region.

In addition to studying the identified region, we also consider the problem of selecting moments and models in the context of (1.1), where only a subset of the moment inequalities are to be used and the true parameter vector $\theta_0$ is assumed to follow a submodel allowing only some selected components to be nonzero (which can be, e.g., the regression coefficients of some selected explanatory variables). Andrews and Soars (2007) employ a modified moment selection procedure to determine which moment inequalities are not binding, by minimizing an information-type criterion. The moment/model selection problem we consider here is different from theirs. In this paper, we have two goals in the selection procedure: first, selecting a true moment/model and second, among the true candidates, selecting the "optimal" one, in a sense which will be described in Section 4. Since the true parameter is not point-identified, it is impossible to test the moment inequalities evaluated at the true parameter. Hence, the moment inequalities in this paper are true in the sense that, fixing the dimension of the parameter vector and the parameter space, the identified region defined by these moment inequalities on the fixed parameter space is *not* empty. In addition, we observe that whether a set of moment inequalities is satisfied or not also depends on the parameters that are included in the model and hence is related to the parameter space. In some situations, a set of moment inequalities defines a valid (nonempty) identified region on one parameter space, but cannot if one or more of the parameters are excluded from the model. This is a result of the reduction of the dimension of the parameter space. By treating the set of moments and the set of nonzero parameters as a combination, the problems of moment selection and model selection are combined. In Section 4, we propose the maximum posterior criterion (MPC) to select the combination that has the largest posterior probability.

We are interested in examining the asymptotic properties of the model/moment combinations proposed by the maximum posterior. We hope that the maximum posterior criterion will produce a desirable combination in the following three senses. First, asymptotically, it should be true. Second, it is desirable that it should contain as many moment inequalities as possible since, intuitively, the more moment inequalities we have, the smaller the identified region is and hence we have more information about the true parameter. Finally, the model should be as simple as possible, that is, the parameter subspace should have the smallest dimension. We show in Section 4 that, indeed, with suitable specifications on the prior, the maximum posterior criterion can produce such a desirable combination with probability tending to one as the sample size increases to infinity. Such a result will be referred to as the consistency of the MPC for model/moment selection.

The remainder of this paper is organized as follows. Section 2 describes a general moment inequality model and the construction of the limited information likelihood. We also provide a general consistency theorem on set estimation based on the posterior c.d.f. function. Section 3 provides a detailed large-sample analysis of the behavior of the posterior distribution. In particular, we will derive the convergence rates both within the interior of $\Omega$ and on any $\delta$-contraction outside $\Omega$. Section 4 studies the problem of moment/model selection. Section 5 displays some simulation results. Finally, Section 6 concludes with a discussion. Proofs are given in Appendices A–C.

## 2. Moment inequality models.

2.1. *Limited information likelihood.* Suppose that for $\theta \in \mathbb{R}^d$, we have moment inequality conditions

$$Em_j(X_i, \theta) \geq 0, \qquad j = 1, \ldots, p.$$

Let

$$m(X, \theta) = (m_1(X, \theta), m_2(X, \theta), \ldots, m_p(X, \theta))^T.$$

The moment inequalities can then be rewritten as

(2.1) $\qquad Em(X, \theta) = \lambda \qquad$ for some $\lambda \in [0, \infty)^p$.

Here, $\theta$ is the structural parameter of interest, for example, $\theta = EY$, the mean of the unobserved random variable $Y$ in Examples 1.1 and 1.2, and $\lambda$ is the bias parameter of $Em(X, \theta)$, for example, $\lambda = (EY_2 - \theta, \theta - EY_1)^T$ in Example 1.1. Let $\theta_0$ be the true parameter value of $\theta$ and $\lambda_0$ be the true bias parameter when $\theta = \theta_0$. Suppose that the prior of $\theta_0$ is supported on a large enough compact set that contains the identified region. We are interested in constructing the marginal posterior for $\theta_0$.

In addition, let $\bar{m}(\theta) = \frac{1}{n}\sum_{i=1}^{n} m(X_i, \theta)$ and $G(\theta, \lambda) = \bar{m}(\theta) - \lambda$. Then, after the bias parameter $\lambda$ is introduced, $G$ can be considered as the "debiased" sample moment. In other words, $G$ is an estimating function with $EG(\theta, \lambda) = 0$. It is overparametrized, meaning that the dimension of $(\theta, \lambda)$ is greater than the dimension of $G$, and, hence, we cannot consistently estimate $\theta_0$ by solving $G(\theta, \lambda) = 0$ directly.

Under some regularity conditions, by the central limit theorem,

(2.2) $\qquad\qquad \sqrt{n}G(\theta, \lambda)|_{\theta=\theta_0, \lambda=\lambda_0} \xrightarrow{d} N_p(0, V_0),$

where $V_0 = \mathrm{Var}(m(X,\theta_0))$. We can therefore formally construct a "likelihood" function:

$$(2.3) \qquad p(X^n|\theta,\lambda) = \frac{1}{\sqrt{\det(2\pi/nV_0)}} e^{-n/2 G(\theta,\lambda)^T V_0^{-1} G(\theta,\lambda)}.$$

Note that for $\theta \neq \theta_0$, (2.2) is not true in general. In fact, we cannot find a $\lambda \in [0,\infty)^p$ such that $Em(X,\theta) = \lambda$ for $\theta \notin \Omega$. Hence, (2.3) is not the large-sample conditional p.d.f. of $G$ for general $(\theta,\lambda)$. The asymptotic result (2.2) alone would not allow us to derive a likelihood function over the entire $\Theta \times [0,\infty)^p$. To solve this problem, Kim (2002) introduced the concept of *limited information likelihood*. For each parameter $\theta \in \Theta$, although (2.3) may not be the true probability density of $X^n$, it is shown to be proportional to the density that is closest to the true density in the Kullback–Leibler distance, among a family of densities satisfying the moment condition $EG(\theta,\lambda) = 0$. The "likelihood" in (2.3) is therefore the limited information likelihood of $\theta$ and $\lambda$, which is the best approximation to the true density that satisfies the moment restrictions. The concept of the Kullback–Leibler information distance and applications of it can be found in a number of works, such as Cover and Thomas (1991) and Zellner (1994).

Let $p(\lambda)$ be the marginal prior of $\lambda$. Assume that $\lambda$ and $\theta$ are independent, that is, the conditional prior of $\lambda$ given $\theta$ is equal to the marginal prior of $\lambda$. Since we are only interested in $\theta$, we thus integrate out $\lambda$ to obtain the limited information likelihood function for $\theta$:

$$(2.4) \qquad \begin{aligned} L(\theta) &= p(X^n|\theta) \\ &= \int_{[0,\infty)^p} p(X^n|\theta,\lambda) p(\lambda|\theta) \, d\lambda \\ &= \int_{[0,\infty)^p} p(X^n|\theta,\lambda) p(\lambda) \, d\lambda. \end{aligned}$$

The fact that $\lambda$ is a location parameter of (2.3) makes the problems solvable. This will be described in detail in Section 3.1.

In practice, the asymptotic variance $V_0$ in (2.3) is not known, but it can be shown to have very little influence on the inference about $\theta$ in the current situation of partially identified moment inequality models. In future expositions, we will replace $V_0$ by a prespecified nonsingular matrix $V$ and show that $L(\theta)$ has good (and very similar) frequentist properties for inference on $\theta$, whatever $V$ is chosen. (A more delicate treatment would be to approximate $V_0$ by a sample analog and replace the true parameter $\theta_0$ in $V_0$ by the unknown argument $\theta$. This will be left for future work. We expect that similar techniques will lead to similar results in this treatment, but the technical details can be much more complicated.)

2.2. *A general result on the posterior set estimation.* We first define some notation that will be used subsequently. Throughout this paper, let $A^c$ and $\mathrm{int}(A)$ denote the complement and interior of a set $A$, respectively. In addition, following CHT's notation, $\forall \delta > 0$, let $(\Omega^c)^{-\delta}$ be the $\delta$-contraction of $\Omega^c$:

$$(\Omega^c)^{-\delta} = \{\theta \in \Theta : d(\theta, \Omega) \geq \delta\}.$$

Let $B(\omega, r)$ denote an open ball around $\omega$: $B(\omega, r) = \{\theta : d(\omega, \theta) < \delta\}$, where $d(\omega, \theta)$ denotes the Euclidean distance between $\omega, \theta$. Let $d_H(A, B)$ denote the Hausdorff distance between sets $A$ and $B$:

$$d_H(A, B) = \max\Big\{\sup_{a \in A} d(a, B), \sup_{b \in B} d(A, b)\Big\},$$

where $d(a, B) = \inf_{b \in B} d(a, b)$. We say that a set estimator $A_n$ consistently estimates $\Omega$ if

$$d_H(A_n, \Omega) \to 0 \quad \text{in probability.}$$

Moreover, for two sequences $\{a_n\}_{n=1}^\infty$ and $\{b_n\}_{n=1}^\infty$, we write $a_n \succ b_n$ if $\frac{a_n}{b_n} \to \infty$. Finally, we write "w.p.a.1" as shorthand for "with probability approaching one in the probability distribution of $X^n$ as $n \to \infty$."

Let $p(\theta)$ be the prior of $\theta$. By Bayes' rule, the posterior of $\theta$ then satisfies

(2.5) $$p(\theta | X^n) \propto p(\theta) L(\theta).$$

It is desirable for the posterior to possess some "good" frequentist properties. Roughly speaking, we want to see that the posterior density of $\theta$ concentrates near $\Omega$ and drops dramatically to zero outside $\Omega$, with a high probability as $n$ increases. The significant difference of such asymptotic behavior between the inside and outside of the identified region implies that the resulting posterior has the capability to produce consistent set estimation for $\Omega$. Such a relation between a "good" posterior and its capability to estimate $\Omega$ is demonstrated below for a scalar function of $\Omega$. (A more general estimation of $\Omega$ itself will also be discussed in Section 3.)

The posterior probability that $\theta$ belongs to a set $A$ is

$$P(\theta \in A | X^n) = \int_A p(\theta | X^n) \, d\theta.$$

DEFINITION 2.1 (Dense). A subset $A \subset \Omega$ is said to be *dense* in $\Omega$ if $\forall \omega \in \Omega \setminus A$ and any neighborhood $U_w$ of $\omega$, we have $U_w \cap A \neq \phi$.

An equivalent definition of dense subsets in real analysis is that the closure of $A$ is $\Omega$, that is, $\mathrm{cl}(A) = \Omega$. We will consider the large-sample behavior of the posterior distribution on a dense subset of $\Omega$.

Suppose that, instead of $\theta$, we are interested in the functions of $\theta$, $g(\theta)$, where $g:\Theta \to \mathbb{R}$ is some known continuous mapping. For instance, if we are interested in the $i$th component of $\theta$, then $g(\theta) = \theta_i$. Let $g(\Omega) = \{g(\theta): \theta \in \Omega\}$, the image of $g$. We are interested in estimating $g(\Omega)$ directly. We impose the following assumptions.

ASSUMPTION 2.1. $\Theta$ is compact.

ASSUMPTION 2.2. $\Omega$ is compact and connected.

In moment inequality models, the compactness of $\Omega$ follows from assuming $Em_j(X,\cdot):\Theta \to \mathbb{R}$ to be continuous for each $j$. We assume $\Omega$ to be connected so that the intermediate value theorem on a topological space holds.

ASSUMPTION 2.3. $g:\Theta \to \mathbb{R}$ is continuous on $\Theta$.

The estimation of $g(\Omega)$ to be constructed is based on the inverted posterior c.d.f. of $g(\theta)$. Let $F_g(x) = P(g(\theta) \leq x|X^n)$, the posterior c.d.f. of $g(\theta)$. Let

$$F_g^{-1}(y) = \inf\{x: F_g(x) \geq y\}.$$

Then $x \geq F_g^{-1}(y)$ if and only if $F_g(x) \geq y$. The following theorem provides a general consistency result of a set estimator of $g(\Omega)$ based on the posterior c.d.f. Note that since it can be shown that $g(\Omega) = [\inf_{\theta \in \Omega} g(\theta), \sup_{\theta \in \Omega} g(\theta)]$, one might think that a more natural set estimator can be constructed by finding estimators for the end points of the interval $g(\Omega)$. This idea works, for example, in Example 1.1, where $[EY_1, EY_2]$ can be estimated by $[\bar{Y}_1, \bar{Y}_2]$, as both $Y_1$ and $Y_2$ are observable. However, in a more general setting, estimating the end points $\inf_{\theta \in \Omega} g(\theta)$ and $\sup_{\theta \in \Omega} g(\theta)$ would require estimating $\Omega$ first. The estimator proposed in the following theorem provides a way of estimating the interval directly.

THEOREM 2.1. *Under Assumptions 2.1–2.3, assume that there exists* $\{\pi_n\}_{n=1}^\infty$, $\pi_n \to 0$, *such that:*

1. $\forall \delta > 0$, $P(\theta \in (\Omega^c)^{-\delta}|X^n) = o_p(\pi_n)$;
2. *there exists a dense subset $A \subset \Omega$ such that $\forall \omega \in A$, $\exists R_w > 0$ such that when $\rho < R_w$, $P(\theta \in B(\omega, \rho)|X^n) \succ \pi_n$ w.p.a.1.*

*If we let $\hat{g} = [F_g^{-1}(\pi_n), F_g^{-1}(1-\pi_n)]$, then*

$$d_H(\hat{g}, g(\Omega)) \to 0 \quad \text{in probability.}$$

REMARK 2.1.

1. The consistent set estimator depends on the choice of $\pi_n$. However, we do not pursue an operational way of constructing the estimator based on the posterior distribution in this paper because there are many frequentist methods to achieve this purpose, for instance, CHT, Beresteanu and Molinari (2008), etc. This paper is more focused on the posterior distribution itself. The purpose of this theorem is to demonstrate that the posterior can be used to consistently estimate the identified region, if needed. The posterior distribution can actually provide more information than the identified region when taking into account the prior.
2. We can also provide an exact credible region (based on, say, setting $\pi_n = 0.025$, for instance) for the true parameter, conditional on the observed data. This is parallel to the provision of the confidence intervals with required coverage probabilities in the frequentist approaches of Imbens and Manski (2004), Rosen (2008), etc.
3. It is possible to get an optimal rate of $\pi_n$ for optimal convergence rate in Hausdorff distance. We leave this for future work.

In the next section, we will see that, under some regularity conditions, the posterior distribution of $\theta$ satisfies conditions 1 and 2 of this theorem, which describe the frequentist properties of the posterior. In addition, we will also propose a consistent estimator for $\Omega$, directly based on the log-posterior density.

**3. Posterior properties of moment inequality models.** In this section, we assume that the identified region contains a nonempty interior $\text{int}(\Omega)$. Assuming it is dense in $\Omega$, it is then of interest to study the asymptotic properties of the posterior distribution inside $\text{int}(\Omega)$.

3.1. *The posterior density.* Following the discussions in Section 2.1, we will study a limited information likelihood for $\theta$ defined by

$$(3.1) \quad L(\theta) = \int_{[0,\infty)^p} \frac{1}{\sqrt{\det(2\pi V/n)}} e^{-n/2(\bar{m}(\theta)-\lambda)^T V^{-1}(\bar{m}(\theta)-\lambda)} p(\lambda)\, d\lambda,$$

where $V$ is some preselected positive definite matrix that does not depend on $\theta$. We will use a multivariate exponential distribution as the prior on $\lambda$ throughout this paper:

$$p(\lambda) = \left(\prod_{i=1}^{p} \psi_i\right) e^{-\psi^T \lambda}, \qquad \psi = (\psi_1, \ldots, \psi_p)^T \in [0,\infty)^p, \lambda \in [0,\infty)^p,$$

where $\psi$ is prespecified. We use the exponential prior for ease of integration over $\lambda$. More general choices of $p(\lambda)$ may not allow the integration to be carried out analytically, but the large-sample behavior of the posterior should remain unchanged.

Let $Z_\theta$ be a $p$-dimensional multivariate normal random vector with mean $(\bar{m}(\theta) - \frac{V\psi}{n})$ and variance-covariance matrix $\frac{V}{n}$. A straightforward calculation of (3.1) then leads to

$$(3.2) \qquad L(\theta) = P(Z_\theta \geq 0) e^{-\psi^T \bar{m}(\theta) + 1/(2n)\psi^T V \psi} \left( \prod_{i=1}^{p} \psi_i \right)$$

and we have $p(\theta|X^n) \propto p(\theta)L(\theta)$.

For large values of $n$, by the uniform weak law of large numbers, $\bar{m}(\theta)$ is bounded on $\Theta$ w.p.a.1. Thus, for fixed $\psi$ and $V$, $e^{-\psi^T \bar{m}(\theta) + 1/(2n)\psi^T V \psi}(\prod_{i=1}^{p} \psi_i)$ is bounded away from zero and infinity. Therefore, the only term that characterizes the large-sample properties of the posterior should be $P(Z_\theta \geq 0)$. Moreover, the variance-covariance matrix of $Z$ has order $O_p(n^{-1})$, so we would expect that $\lim_{n \to 0} P(Z_\theta \geq 0) = 1$ in probability if and only if $\bar{m}(\theta) - \frac{V\psi}{n} \geq 0$ w.p.a.1. This depends on whether or not $\theta$ belongs to $\Omega$. For large $n$, the posterior density is positive inside $\Omega$ and drops to zero exponentially fast as $\theta$ goes away from $\Omega$. We will formally examine these asymptotic properties in the next section and will also derive the convergence rate of the posterior probabilities.

3.2. *Large-sample analysis.* We now provide a large-sample analysis of the posterior distribution of the parameter $\theta$.

ASSUMPTION 3.1. $\text{int}(\Omega)$ is nonempty and is dense in $\Omega$.

The assumption that $\text{int}(\Omega)$ is dense in $\Omega$ can be restated as follows: for any $\omega$ on the boundary of $\Omega$ and any neighborhood $U_w$ of $\omega$, $U_w$ contains points in $\text{int}(\Omega)$. Most of the identified regions characterized by moment inequalities possess such a property. We will comment on the case when $\text{int}(\Omega)$ is empty in the discussion section.

ASSUMPTION 3.2. $Em_j(X, \cdot): \Theta \to \mathbb{R}$ is continuous for each $j = 1, \ldots, p$.

This assumption guarantees that $Em(X, \theta)$ is bounded in any compact set and that the uniform law of large number holds. The next assumption puts a regularity condition on the prior of $\theta$.

ASSUMPTION 3.3. $p(\theta)$ is continuous and bounded away from zero and infinity on $\Omega$.

Let $v_{jj}$ be the $j$th diagonal element of $V$. We can write

$$\Omega^c = \left\{ \theta : \min_j Em_j(X, \theta) < 0 \right\} = \left\{ \theta : \min_j \frac{Em_j(X, \theta)}{\sqrt{v_{jj}}} < 0 \right\}.$$

For any $\delta > 0$, let

$$A_\delta = \left\{\theta : \min_j \frac{Em_j(X,\theta)}{\sqrt{v_{jj}}} < -\delta\right\}.$$

Apparently, $A_\delta \subset \Omega^c$.

LEMMA 3.1. *Under Assumptions 2.1, 2.2 and 3.2, if there exists some $a_n \to 0$ such that $\forall \delta > 0$, $P(\theta \in A_\delta | X^n) = o_p(a_n)$, then $\forall \varepsilon > 0$, $P(\theta \in (\Omega^c)^{-\varepsilon} | X^n) = o_p(a_n)$.*

THEOREM 3.1. *Under Assumptions 2.1, 2.2 and 3.1–3.3:*

1. *$\forall \delta > 0$, for some $\alpha > 0$,*

$$P(\theta \in (\Omega^c)^{-\delta} | X^n) = o_p(e^{-\alpha n});$$

2. *$\forall$ nonempty open sets $\Xi \subset \Omega$, in probability,*

$$\liminf_{n \to \infty} P(\theta \in \Xi | X^n) > 0.$$

Hence, we are able to distinguish the asymptotic behavior of the posterior: for large values of $n$, the posterior density is only supported on a neighborhood of the identified region and the posterior distribution drops to zero exponentially fast on any subset that is separated from $\Omega$. Based on these findings, we can construct consistent estimators for both $\Omega$ and its continuous mappings. For the latter task, we can now apply Theorem 2.1. Suppose that $g(\cdot)$ is a continuous real-valued function on $\Theta$ and let $F_g^{-1}(y)$ be the $y$-quantile of the posterior c.d.f. of $g(\theta)$.

THEOREM 3.2. *Under Assumptions 2.1–2.3 and 3.1–3.3, for any sequence $\pi_n = o_p(1)$ such that $\forall a > 0$, $e^{-an}/\pi_n \to 0$, we have*

$$d_H([F_g^{-1}(\pi_n), F_g^{-1}(1-\pi_n)], g(\Omega)) \to 0 \qquad \text{in probability.}$$

We can also consistently estimate $\Omega$ directly using the posterior density function. The consistency is based on the fact that the posterior density attains its peak inside $\Omega$ and is asymptotically supported on the entire identified region. In addition, it drops to zero outside $\Omega$ at an exponential rate. Therefore, by properly choosing a cut-off value $\varepsilon_n$, the region where the log-posterior p.d.f. exceeds its peak minus $\varepsilon_n$ should eventually converge to $\Omega$.

THEOREM 3.3. *Under Assumptions 2.1–2.3 and 3.1–3.3, let $1 \prec \varepsilon_n \prec n$. If we define*

$$A_n = \left\{\theta : \max_{\omega \in \Theta} \ln p(\omega | X^n) - \ln p(\theta | X^n) \leq \varepsilon_n\right\},$$

*then*

$$d_H(A_n, \Omega) \to 0 \qquad \text{in probability.}$$

REMARK 3.1. The estimation established in Theorem 3.3 is easy to implement for the following reasons.

1. Note that

$$\max_{\omega \in \Theta} \ln p(\omega|X^n) - \ln p(\theta|X^n)$$

$$= \max_{\omega \in \Theta} \left( \ln p(\omega) L(\omega) - \ln \int_\Theta p(\theta) L(\theta) \, d\theta \right)$$

$$- \left( \ln p(\theta) L(\theta) - \ln \int_\Theta p(\theta) L(\theta) \, d\theta \right)$$

$$= \max_{\omega \in \Theta} \ln p(\omega) L(\omega) - \ln p(\theta) L(\theta).$$

Therefore, there is no need to normalize $p(\theta)L(\theta)$, avoiding numerical integration of $p(\theta)L(\theta)$.
2. Maximizing $\ln p(\theta)L(\theta)$ is computationally workable since the maxima is attained only inside $\Omega$, where $p(\theta)L(\theta)$ is quite smooth. Hence, Newton–Raphson's algorithm can carry out the maximization.
3. If we set $a_n = \max_{\omega \in \Theta} \ln p(\omega|X^n) - \varepsilon_n$, then $A_n = \{\theta : \ln p(\theta) \geq a_n\}$. The boundary $\{\theta : \ln p(\theta) - a_n = 0\}$ is a closed curve with dimension $d - 1$.

**4. Moment and model selection.** In this section, we discuss the problem of moment and model selection. Suppose that we have $p$ candidate moment inequalities

$$Em_j(X, \theta) \geq 0, \qquad j = 1, \ldots, p,$$

with a $k$-dimensional parameter vector $\theta = (\theta_1, \ldots, \theta_k)^T \in \Theta_1 \times \cdots \times \Theta_k$. The moment selection problem refers to selecting the best subset of the moment inequalities among all of the possible candidates (where there is some notion of optimality), while the model selection procedure addresses the problem of selecting the best model among all subsets of the parameter space where some components of the parameter are set to zero. The possible moment inequalities and corresponding subsets of the parameter space are known. What is not known is which ones are the best.

Instead of selecting the moment inequalities and the parameter subspace as two separate procedures, we select them simultaneously, as a combination. However, there are still two problems to consider: selecting the true combination and, among the true combinations, selecting the optimal one, in the sense that it should contain as many moment inequalities and as few structural parameters as possible. The selection procedure is based on the maximum posterior criterion (MPC): we assign prior probabilities to each candidate moment/model and then derive the posterior probabilities based on the limited information likelihood described in Section 2, by integrating out the structural and nuisance parameters $(\theta, \lambda)$. Finally, the optimal combination is the one with the largest maximum posterior probability. We will examine the asymptotic property of the optimal combination by establishing the consistency of MPC. By consistency, we mean w.p.a.1, MPC will select the true combination with the most moment inequalities and fewest structural parameters (the simplest parameter subspace).

4.1. *Selecting the true combination.* Because of the feature of partial identification, it is impossible to test whether each candidate moment inequality is true at the true parameter. Given a set of moment inequalities and a parameter space, we can tell whether the moment inequalities define a nonempty identified region on the parameter space.

EXAMPLE 4.1 (Interval regression model). Suppose that an interval regression model provides moment inequalities as follows:

$$E(Z_1 Y_1) \leq E(Z_1 Y) \leq E(Z_1 Y_2);$$
$$E(Z_2 Y_1) \leq E(Z_2 Y) \leq E(Z_2 Y_2).$$

We assume that the data-generating process is $E(Y|X) = X^T \theta_0$, where $\theta_0 = (0.9, 0)^T$ is the true parameter and $X_1 \sim \text{uniform}[-1, 1]$, $X_2 = 1$ a.s. Furthermore, let $\binom{Z_1}{Z_2} = \binom{X_1+1}{1}$ and $\binom{Y_1}{Y_2} = \binom{Y+0.1(U_1-1)}{Y+0.1(U_2+1)}$, where $U_1$ and $U_2$ are uniform $[-1, 1]$ independently. $\theta_0$ is not partially identified by the moment inequality models. If we let $(\theta_1, \theta_2)^T \in \mathbb{R}^2$ be the parameter vector, then we have four moment inequalities:

(4.1) $$\tfrac{1}{3}\theta_1 + \theta_2 \geq 0.2;$$

(4.2) $$\tfrac{1}{3}\theta_1 + \theta_2 \leq 0.4;$$

(4.3) $$\theta_2 \geq -0.1;$$

(4.4) $$\theta_2 \leq 0.1.$$

The region defined by (4.1)–(4.4) on $\mathbb{R}^2$ gives the nonempty identified region for $\theta_0$ with a parallelogram boundary. However, if we set $\theta_1 = 0$, (4.1) contradicts (4.4). Hence, in this case, (4.1)–(4.4) defines an empty region.

Let us define a combination $C_s = (M_{s_1}, \Theta_{s_2})$ with a vector index $s = (s_1, s_2)$, $s_1 \in \{1, 2, \ldots, 2^p - 1\}$ and $s_2 \in \{1, \ldots, 2^k\}$. Here, $M_{s_1}$ denotes a subset of moments, for instance, $M_{s_1} = \{m_1\}$, $M_{s_1} = \{m_1, m_2\}$, etc. There are then $2^p - 1$ such possible subsets. In addition, we denote by $\Theta_{s_2}$ the parameter subspace corresponding to the selected model. By definition, $\Theta_{s_2}$ is the subset of vectors with one or more components fixed to be zero. There are $2^k$ possible $\Theta_{s_2}$'s. (Note that we can select none of the parameters, in which case the model is a reduced model, e.g., in Cox's proportional hazard model; if all of the parameters are set to be zero, we get the baseline model.) The combination $C_s$ combines both the candidate moment functions and the parameter subspace. When selecting a subset of moment inequalities, we also specify a subspace of the structural parameter.

EXAMPLE 4.2 (Example 4.1 continued). Let $\Theta_1 \times \Theta_2$ be the parameter space for $(\theta_1, \theta_2)$, chosen large enough so that $\{(\theta_1, \theta_2) : 0.2 \leq \frac{1}{3}\theta_1 + \theta_2 \leq 0.4, -0.1 \leq \theta_2 \leq 0.1\} \subset \Theta_1 \times \Theta_2$. A scope of candidate combinations can be any of the following:

$\{E(Z_1 X^T \theta - Z_1 Y_1)\}, \qquad \Theta_1 \times \Theta_2;$

$\{E(Z_1 X^T \theta - Z_1 Y_1), E(Z_1 Y_2 - Z_1 X^T \theta)\}, \qquad \Theta_1 \times \Theta_2;$

$\{E(Z_2 X^T \theta - Z_2 Y_1)\}, \qquad \{0\} \times \Theta_2;$

$\{E(Z_1 X^T \theta - Z_1 Y_1), E(Z_1 Y_2 - Z_1 X^T \theta), E(Z_2 Y_2 - Z_2 X^T \theta)\}, \qquad \Theta_1 \times \{0\};$

$\vdots$

$\{E(Z_2 Y_2 - Z_2 X^T \theta)\}, \qquad \Theta_1 \times \Theta_2.$

DEFINITION 4.1. A combination $C_s = (M_{s_1}, \Theta_{s_2})$ is true if and only if

$$\inf_{\theta \in \Theta_{s_2}, \lambda \in [0, \infty)^m} \|EM_{s_1}(X, \theta) - \lambda\|^2 = 0,$$

where $m$ denotes the number of candidate moment functions in $M_{s_1}$.

If we let $\Omega_s = \{\theta \in \Theta_{s_2} : EM_{s_1}(X, \theta) \geq 0\}$ be the identified region defined by $C_s$, then this definition is equivalent to saying that $\Omega_s$ is not empty.

We place a discrete prior $p(C_s)$ on all of the candidate combinations. In practice, such a prior can be either uniform [i.e., $p(C_s) = \frac{1}{2^k(2^p - 1)}$ for all $C_s$] or model-dependent, or obtained from previous studies. As in the previous sections, let

$$\lambda = EM_{s_1}(X, \theta)$$

with $\dim(\lambda) = m$ and use the following prior conditional on $C_s$:

$$p(\lambda|C_s) = \prod_{i \leq m} \psi_i e^{-\psi^T \lambda}, \qquad \lambda > 0,$$

where the $\psi_i$'s are the prespecified second-stage parameters. Let $p(\theta|C_s)$ be the conditional prior of the parameter $\theta \in \Theta_{s_2}$. The conditional limited information likelihood is given by

$$L(X^n|\theta,\lambda,C_s) = \frac{1}{\sqrt{\det(2\pi/nV)}} e^{-n/2(\bar{M}_{s_1}(\theta)-\lambda)^T V^{-1}(\bar{M}_{s_1}(\theta)-\lambda)},$$

where $\bar{M}_{s_1}(\theta) = \frac{1}{n}\sum_{i=1}^n M_{s_1}(X_i,\theta)$. The posterior of $C_s$ can then be obtained by integrating out $\theta$ and $\lambda$, which is proportional to the "integrated likelihood,"

$$(4.5) \quad p(C_s|X^n) \propto \iint_{\Theta_{s_2} \times [0,\infty)^m} L(X^n|\theta,\lambda,C_s) p(\theta|C_s) p(\lambda|C_s) p(C_s)\, d\theta\, d\lambda.$$

A remark on the "$d\theta$" part of this integration: the integration is with respect to the nonzero elements of $\theta \in \Theta_{s_2}$, where $\Theta_{s_2}$ is the parameter space of those free parameters only. For instance, suppose that the full parameter is $(\theta_1,\theta_2,\theta_3) \in \Theta_1 \times \Theta_2 \times \Theta_3$. Once we set $\theta_3 = 0$, then $\Theta_{s_2} = \Theta_1 \times \Theta_2$ and integrating over $\theta$ becomes a two-dimensional integration (w.r.t. $\theta_1,\theta_2$). Otherwise, if we set $\Theta_{s_2} = \Theta_1 \times \Theta_2 \times \{0\}$ and still treat it as a three-dimensional integration, $\Theta_{s_2}$ would have a zero Lebesgue measure and, as a result, the integration would always be zero.

We make the following assumptions.

ASSUMPTION 4.1. The parameter space $\Theta$ of the full model is compact.

The next assumption is imposed on the prior of $\theta$.

ASSUMPTION 4.2. If $C_s$ is true, then $p(\theta \in \Omega_s|C_s) > 0$.

ASSUMPTION 4.3. $p(C_s) > 0$ for each combination $C_s$.

The following assumption implies that $\Omega_s$ is nonempty and compact, given that $C_s$ is true.

ASSUMPTION 4.4. $Em_j(X,\theta)$ is continuous on $\Theta_s$ for each $m_j$ in $C_s$.

Intuitively, we should select as many moment inequalities as possible since the more moment inequalities there are, the smaller the identified region is. However, if one or more of the selected moment inequalities are false, the identified region is empty. The following theorem illustrates that the posterior probability is exponentially small if the selected combination is false.

THEOREM 4.1. *For combination $C_s = (M_{s_1}, \Theta_{s_2})$, under Assumptions 4.1–4.4:*

1. *if $C_s$ is true, then, in probability*
$$\liminf_{n \to \infty} p(C_s | X^n) > 0;$$

2. *if $C_s$ is not true, then for some $\alpha > 0$,*
$$p(C_s | X^n) = o_p(e^{-\alpha n}) p(C_s).$$

4.2. *Selecting the optimal combination.* The maximum posterior procedure provides an optimality criterion to select the combination with the largest posterior probability

(4.6) $$C^* = \arg\max_{C_s} p(C_s | X^n).$$

We are interested in studying the asymptotic properties of the optimal $C^*$. We hope that the MPC will produce a desirable combination in the following three senses. First, asymptotically, $C^*$ should be true. Second, it is desirable that it should contain as many moment inequalities as possible since, intuitively, the more moment inequalities we have, the smaller the identified region is and hence we have more information about the true parameter. Finally, the model should be as simple as possible.

We impose the following assumption in addition to Assumptions 4.1–4.4.

ASSUMPTION 4.5. Each true candidate combination $C_s$ defines a connected $\Omega_s$.

We first consider using a (discrete) uniform prior for the candidate combinations: for all $C_s$,

(4.7) $$p(C_s) = \frac{1}{2^k(2^p - 1)}.$$

Although this seems to be a natural prior to use, we found examples where it actually functions undesirably for model or moment selection. For example, suppose that we want to compare the posterior probabilities of two candidate combinations, $C_s^1 = (M_{s_1}^1, \Theta_{s_2}^1)$ and $C_s^2 = (M_{s_1}^2, \Theta_{s_2}^2)$, using the Bayes factor
$$\mathrm{BF}_{12} = \frac{p(C_s^1 | X^n)}{p(C_s^2 | X^n)}.$$
We fix $\Theta_{s_2}^1 = \Theta_{s_2}^2$ and assume $M_{s_1}^1 \subset M_{s_1}^2$, that is, the moment inequalities of $M_{s_1}^1$ are strictly contained in the moment inequalities of $M_{s_1}^2$. If both $C_s^i$, $i = 1, 2$, are true, then the identified region $\Omega_s^2$ defined by $C_s^2$ should be a strict subset of $\Omega_s^1$ defined by $C_s^1$. As explained before, a smaller identified region is preferable since it provides more precise information about the true parameter. Hence, we hope that $\text{BF}_{12}$ is asymptotically less than one if the MPC criterion is consistent with this intuition. However, the following theorem indicates that, if the uniform prior (4.7) is specified, then the result is in exactly the opposite direction.

For a matrix $M$, define $\|H\| = \sqrt{\text{trace}(HH^T)}$.

THEOREM 4.2. *Suppose that Assumptions 4.1–4.5 are satisfied and a uniform prior (4.7) is applied. Suppose that both $C_s^i$, $i = 1, 2$, are true, $M_{s_1}^1 \subset M_{s_1}^2$ and $\Theta_{s_2}^1 = \Theta_{s_2}^2$, In addition, suppose that $\psi_i$ satisfies $\psi_i < e^{-\|\psi\| \cdot \sup_\Theta \|Em(X, \theta)\|}$, $i = 1, \ldots, p$. Then, w.p.a.1,*

$$\text{BF}_{12} > 1.$$

Here, $\psi_i$ is the second-stage parameter of the exponential prior of $\lambda_i$. In practice, a small $\psi$ is preferable because it leads to a noninformative prior on $\lambda$. However, Theorem 4.2 says that if $\psi$ is small (satisfying $\psi_i < e^{-\|\psi\| \cdot \sup_\Theta \|Em(X, \theta)\|}$, $i = 1, \ldots, p$) and a uniform prior for the candidate combinations is used, then the result will be negative: a candidate with fewer moment constraints has a larger posterior distribution. However, this is not a warning about the method, but rather about the potential danger of the seemingly innocent choice of the uniform prior on the candidate combinations. With this prior, it will be shown in Appendix C that the posterior of each true combination is of order $O_p(1)$ and, up to the leading order, is proportional to the prior measure of the identified region, as well as the product of the $\psi_i$'s. As more moment inequalities are added, the identified region gets smaller. Also, more small $\psi_i$'s are added in the product term. Both of these factors make the resulting posterior probability smaller.

We will address this problem either by using a more informative prior on the candidate combinations (Approach 1) or by placing some uninformative priors on some components of the parameters $\theta$ and $\lambda$ (Approach 2). Either way, the posterior of each candidate is no longer of order $O_p(1)$ and the order of the optimal candidate's posterior will be the largest, which overrides the effects from the prior measure of the identified region and the product of the $\psi_i$'s.

4.3. *Prior selection for consistency of MPC.* We propose two approaches to address the problem illustrated in Theorem 4.2.

APPROACH 1. One approach is to change the priors of all the candidate combinations. Instead of the uniform (equally likely) priors, we use unequal priors. Specifically, the priors are data-size-dependent and tend to encourage those combinations with more moment inequalities and simpler parameter structures. One such prior could be

$$p(C_s) \propto n^{\alpha[\dim(M_s) - \dim(\Theta_s)]} \tag{4.8}$$

for some $\alpha > 0$. This choice of prior encourages $C_s$ with large $\dim(M_s) - \dim(\Theta_s)$.

One needs be aware that although $p(C_s) \propto e^{n(\dim(M_s) - \dim(\Theta_s))}$ also reward large values of $\dim(M_s) - \dim(\Theta_s)$, we do not recommend its use. This is because we have shown earlier that $p(C_s|X^n) = o_p(e^{-\alpha n})p(C_s)$, that is, the posterior probability of a false combination is exponentially small. However, if an exponentially large prior is used, it may override the "big gap" between the false and true combination posteriors.

The drawback of the unequal prior (4.8) is that it is not a uniform one. In Bayesian analysis, it is usually the nature of the data that determines the properties of the posterior and the priors are usually chosen to be uninformative. One may consider using another approach to deriving the priors.

APPROACH 2. In this approach, we still use the discrete uniform (equally likely) prior for the candidate combinations. However, we partition the parameters $\theta$ and $\lambda$ into "restricted" and "unrestricted" parts, according to the biases of the selected and unselected moment functions. Formally, let

$$\lambda = EM(X, \theta),$$

where $M(X, \theta) = (m_1(X, \theta), \ldots, m_p(X, \theta))^T$, the vector of all the candidate moments, and $\theta = (\theta_1, \ldots, \theta_k)^T$, the vector of full parameters supported on $\Theta_1 \times \cdots \times \Theta_k$. Suppose that a combination $C_s = (M_{s_1}, \Theta_{s_2})$ selects $m$ moment conditions $M_{s_1}$ and leaves the rest of the moments (denoted by $M_{s_1}^c$) unused, while selecting a submodel parameterized by $\theta_s \in \Theta_{s_2}$, setting all of the other components of $\theta$ (denoted by $\theta_s^c$) to be zero. One can view model selection as placing a restriction on $\theta$, while moment selection can be viewed as placing a restriction on the bias $\lambda$.

Let $\lambda_s$ be the subvector of $\lambda$ corresponding to the selected moments. Let $\lambda_s^c$ denote the remaining components of $\lambda$ corresponding to $M_{s_1}^c$. We then have

$$EM_{s_1}(X, \theta_s) = \lambda_s, \qquad \lambda_s \geq 0,$$
$$EM_{s_1}^c(X, \theta_s) = \lambda_s^c, \qquad \lambda_s^c \in \mathbb{R}^{p-m}.$$

The bias $\lambda_s$ for the selected moments is restricted to be nonnegative, while the bias $\lambda_s^c$ for the unselected moments is left unrestricted. Therefore, we can define *restricted parameters* as $(\lambda_s, \theta_s^c)$, with restrictions

$$\lambda_s \geq 0, \qquad \theta_s^c = 0.$$

In addition, we call the remaining parameters $(\lambda_s^c, \theta_s)$ *unrestricted parameters* because $(\lambda_s^c, \theta_s) \in \mathbb{R}^{p-m} \times \Theta_{s_2}$. [We have thus partitioned the moment functions into $M(X, \theta_s) = (M_{s_1}(X, \theta_s)^T, M_{s_1}^c(X, \theta_s)^T)^T$ and parameters into $\lambda = (\lambda_s, \lambda_s^c)$ and $\theta = (\theta_s, \theta_s^c)$.]

For the unrestricted (selected) parameter $\theta_s$, let $t = \dim(\theta_s)$. We release the compactness assumption on the support of $\theta_s$ and assume it is supported on $\mathbb{R}^t$. We then place the following "working" priors on the unrestricted parameters:

$$(4.9) \qquad p(\lambda_s^c | C_s) \sim N_{p-m}(0, \sigma_n^2 I_{p-m}),$$

$$(4.10) \qquad p(\theta_s | C_s) \sim N_t(0, n\sigma_n^2 I_t),$$

where $N_t$ denotes $t$-dimensional multivariate normal distribution. We require that $\sigma_n \to \infty$ as $n$ tends to infinity, but $\sigma_n / e^{\alpha n} \to 0, \forall \alpha > 0$. Here, $I_t$ denotes the $t \times t$ identity matrix. Since the variance of each component of $\lambda_s^c$ and $\theta_s$ approaches infinity as the sample size tends to infinity, (4.9) and (4.10) tend to be very flat. Hence, this choice of prior is uninformative. In addition, we still assign an exponential prior to the restricted parameter $\lambda_s$,

We then include both selected $M_s$ and unselected $M_s^c$ to construct the limited information likelihood, which depends only on the unrestricted $\theta_s$ since $\theta_s^c = 0$:

$$L(X^n | \theta_s, \lambda, C_s) = \frac{1}{\sqrt{\det(2\pi/nV)}} e^{-n/2(\bar{M}(\theta_s) - \lambda)^T V^{-1}(\bar{M}(\theta_s) - \lambda)},$$

where $\bar{M}(\theta_s) = \frac{1}{n} \sum_{i=1}^n M(X_i, \theta_s)$. The posterior of $C_s$ can then be obtained by integrating out $\theta_s$ and $\lambda = (\lambda_s^T, \lambda_s^{cT})^T$, which is proportional to the "integrated likelihood":

$$p(C_s | X^n) \propto \iint_{\Theta_{s_2} \times [0,\infty)^m \times \mathbb{R}^{p-m}} L(X^n | \theta_s, \lambda, C_s) p(\theta_s | C_s) p(\lambda_s | C_s) p(\lambda_s^c | C_s)$$

$$\times p(C_s) \, d\theta_s \, d\lambda_s \, d\lambda_s^c.$$

Note that since multivariate normal priors are placed on the unrestricted parameters, the parameter vector $\theta_s$ is no longer supported on a compact set. As a result, to derive the large-sample properties of $p(C_s | X^n)$ becomes much harder than in the previous sections because $EM(X, \theta)$ may not be bounded on the noncompact parameter space. Instead of providing a general proof, we will study the problem for a specific model of Example 1.3 because this is the most interesting example in the Introduction where we consider model and moment selection. In this example, the model selection can correspond to selecting the useful explanatory variables and the moment selection can correspond to selecting the valid instrumental variables. The key feature of this example is that the moment inequality functions are linearly dependent on $\theta$. We point out that to establish consistency of Approach 2 in a more general framework is possible, but would require additional assumptions that are much more technical.

ASSUMPTION 4.6. *Suppose that the moment inequalities are given by*

$$EZ(Y_2 - X^T\theta) \geq 0, \qquad EZ(X^T\theta - Y_1) \geq 0,$$

*where $Y_1 \leq Y_2$, and $Z$ is a vector of positive random variables. Assume that:*

 (i) $\operatorname{rank}(EZX^T) = \dim(X)$;
 (ii) *there exists at least one true candidate combination;*
 (iii) *each true candidate $C_s$ defines a compact identified region.*

This assumption rules out those candidates that lead to unbounded identified regions, in which case integrals can be infinite.

The following theorem shows that with either one of the two approaches described above, asymptotically, the optimal $C^*$ can have all of the desirable properties: it is true, it defines the smallest nonempty identified region and it corresponds to the simplest model (with the smallest number of free parameters). We refer to this result as the consistency of maximum posterior criterion for the Bayesian moment/model selection problem.

THEOREM 4.3 (Consistency of MPC). *Let*

$$C^* = \arg\max_{C_s} p(C_s | X^n),$$

*where $p(C_s | X^n)$ is obtained from either one of the following:*

1. *prior (4.8) for candidate combinations, with Assumptions 4.1–4.5;*
2. *flat prior for candidate combinations, and parameter priors (4.9) and (4.10), with Assumptions 4.2–4.6 for the instrumental variable interval regression model (Example 1.3).*

*Then, w.p.a.1, $C^*$ satisfies:*

1. *it is true;*
2. *among all of the true combinations, it has the largest $\dim(M_s) - \dim(\Theta_s)$.*

**5. Monte Carlo experiments.** This section presents some Monte Carlo simulation results. We first provide evidence for the finite-sample behaviors of the consistent estimators described in the previous sections as well as the posterior distribution. We simulate the models described in Examples 1 and 2 in Chernozhukov, Hong and Tamer (2007). We then show simulated evidence of the consistency of MPC for the moment/model selection problem.

TABLE 1
*Estimation based on posterior density*

| $\varepsilon_n$ | $\sqrt{n}$ | $\ln n$ | $\ln \ln n$ |
|---|---|---|---|
| $n = 500$ | $[-0.2841, 5.2634]$ | $[-0.123, 5.113]$ | $[-0.0389, 4.702]$ |
| $n = 1000$ | $[-0.2362, 5.2267]$ | $[-0.1135, 5.0977]$ | $[-0.0342, 4.9110]$ |
| $n = 5000$ | $[-0.1158, 5.1233]$ | $[-0.0477, 5.0476]$ | $[-0.0202, 4.9779]$ |

EXAMPLE 5.1 (Interval data). Consider Example 1.1 described in Section 2. The parameter of interest is $\theta = E(Y)$ with moment inequalities

$$E(Y_2 - \theta) \geq 0, \qquad E(\theta - Y_1) \geq 0.$$

We set $Y_1 \sim N(0, 0.1)$ and $Y_2 \sim N(5, 0.1)$, then $\Omega = [0, 5]$. $Y_1$ and $Y_2$ are generated independently and observations with $Y_1 > Y_2$ are discarded. We also set $\psi_1 = 0.1$, $\psi_2 = 0.5$ and $V = I$, the identity matrix in the likelihood function. In addition, we place a flat prior on $\theta$. We report the estimated identified interval of $\theta$ described both in Theorem 3.2 with $g(\theta) = \theta$ and in Theorem 3.3 for sample sizes $N = 500$, $1000$, $5000$ and various choices of $\varepsilon_n, \pi_n$.

Table 1 reports the estimation of $\Omega$ as in Theorem 3.3. To compare the results corresponding to the choices of $\varepsilon_n$, for each interval $[a, b]$, we calculate $\gamma = (a - 0)^2 + (b - 5)^2$. We find that $\varepsilon = \ln \ln n$ performs better than the other two choices since it has a lower $\gamma$ value.

To construct the estimator based on the posterior c.d.f., we carry out the Metropolis algorithm to draw $B = 5000$ samples from the posterior distribution, then calculate the $\pi_n$-quantile of the empirical c.d.f. with various choices of $\pi_n$. For the Metropolis algorithm, we set initial value $\theta_0 = 1$ and a jump distribution $\tilde{\theta} \sim N(\theta_j, 0.5)$. Table 2 reports the findings with $\pi_n = e^{-\sqrt{n}}, n^{-1}$ and $1/\ln n$. $\pi_n = \frac{1}{n}$ appears to be a better choice compared with other two. We also note that $\pi_n = 1/\ln n$ tends to zero too slow to fully estimate the entire identified interval: the estimated interval shrinks too much inside $\Omega$.

TABLE 2
*Estimation based on empirical c.d.f.*

| $\pi_n$ | $e^{-\sqrt{n}}$ | $\frac{1}{n}$ | $\frac{1}{\ln n}$ |
|---|---|---|---|
| $n = 500$ | $[-0.0716, 5.0418]$ | $[-0.0498, 5.0069]$ | $[0.4048, 3.3447]$ |
| $n = 1000$ | $[-0.0422, 4.9983]$ | $[-0.0383, 5.0164]$ | $[0.3304, 3.2542]$ |
| $n = 5000$ | $[-0.0155, 5.0098]$ | $[-0.0063, 4.9927]$ | $[0.2717, 3.8012]$ |

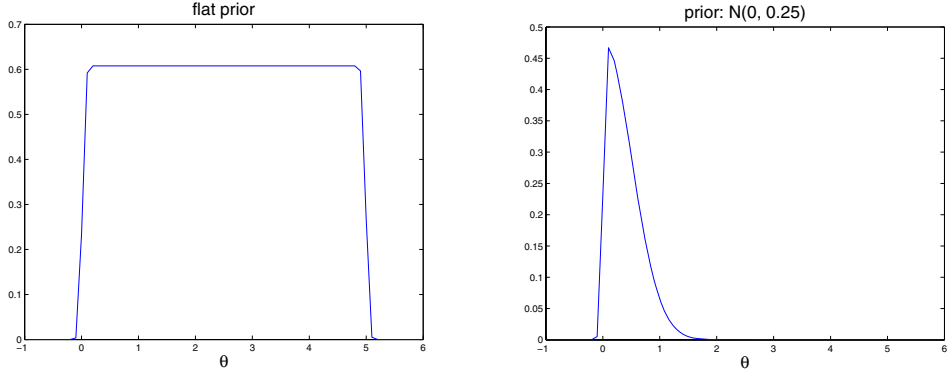

Fig. 1. *Example 5.1: the posterior density function of $\theta$.*

In addition, Figure 1 plots the posterior density function of $\theta$ with two choices of priors: a flat prior and an $N(0, 0.25)$ prior. Theoretically, one needs to truncate the normal distribution so that the priors are supported on a compact set. However, since the tail of the normal density function is very thin and we can choose a very large parameter space, we believe a normal prior is workable here. We see that when a flat prior is used, the posterior density function is high on the entire identified interval $[0, 5]$, but when the prior is set to be $N(0, 0.25)$, most posterior mass falls in $[0, 2]$, which tends to underestimate the true identified interval. However, with this more informative prior, the posterior provides more information about the location of $\theta$.

EXAMPLE 5.2 (Interval outcomes in regression models). We simulate the instrumental inequality model described in Example 1.3,

$$E(ZY_1) \leq E(ZX^T)\theta \leq E(ZY_2),$$

where $\theta = (\theta_1, \theta_2)^T, X = (X_1, X_2)^T$ and $Y = (Y_1, Y_2)^T \in \mathbb{R}^2$. Generate $X \sim N_2((1, 1)^T, I_2)$. Let $Z_1 = X_1 + X_2$ and $Z_2 = X_1 + 2X_2$. Generate $Y_1 \sim N(3, 0.1)$, $Y_2 \sim N(6, 0.1)$ independently. We discard a stack of generated data if either $Z_1$ or $Z_2$ is negative. The identified region is $f\Omega = \{\theta : 2 \leq \theta_1 + \theta_2 \leq 4, 9 \leq 4\theta_1 + 5\theta_2 \leq 18\}$, a two-dimensional region with parallelogram boundary. To estimate this model, set $\psi = (0.1, 0.1, 0.5, 0.5)^T$, $V = I$. Fixing sample size $n = 500$, we implement the Metropolis algorithm to draw $B = 5000$ samples from the posterior distribution.

We first put a flat prior on $\theta$. Figure 2 (left) displays the parallelogram boundary of $\Omega$, as well as 5000 draws from the posterior distribution. Most of the draws fall uniformly inside the identified set, except for those close to the two opposite angles of the parallelogram. We can see that there is a small "bias" at the boundaries.

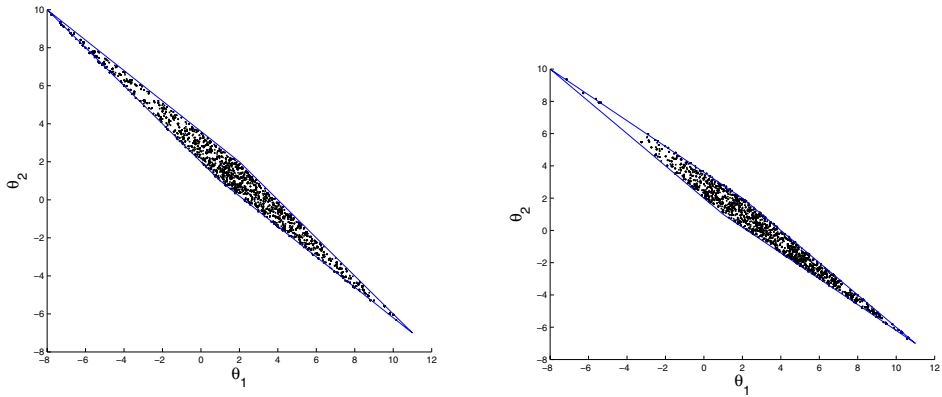

Fig. 2. *Example 5.2: the identified set and MCMC draws. Left: flat prior; right: prior (5.1).*

In order to show that when a more informative prior is applied, the posterior distribution indeed provides more information about the location of the true parameter inside the identified region, we repeat the same MCMC procedure, but with prior distribution

(5.1) $$\theta_1 \sim N(10, 12^2), \qquad \theta_2 \sim N(-6, 12^2),$$

where $\theta_1$ and $\theta_2$ are a priori independent. This prior can be used when, for instance, a previous study estimates that $E\theta_1 \approx 10$ and $E\theta_2 \approx -6$, with the same standard deviation, 12. Figure 2 (right) displays 5000 MCMC draws from the posterior derived from prior (5.1). We see that the draws mostly concentrate at the lower-right corner inside the identified region, which is close to $(10, -6)$, showing that our Bayesian approach indeed provides more information on $\theta$ in this case than the frequentist method; the latter would only estimate the identified region and provide a confidence set, but not tell how $\theta$ is distributed inside it.

EXAMPLE 5.3 (Moment selection: interval censored data). Suppose $\theta \in \Theta \subset \mathbb{R}$. We consider four moment conditions:

(5.2) $$EY_1 \geq \theta;$$

(5.3) $$EY_2 \leq \theta;$$

(5.4) $$EY_3 \leq \theta;$$

(5.5) $$EY_4 \geq \theta.$$

If we assume that $\Theta = [0, 10]$ and $EY_1 < 0 < EY_2 < EY_3 < EY_4 < 10$, then for fixed $\Theta$, (5.2) is incorrect. We generate $N$ i.i.d. data from $Y_1 \sim N(-1, 0.1)$, $Y_2 \sim N(1, 0.1)$, $Y_3 \sim N(2, 0.1)$ and $Y_4 \sim N(3, 0.1)$, with $N = 100, 1000$ and 5000. We fix $\theta \in \Theta$ and use the prior described in Section 4.3, Approach 2 to construct the posterior probabilities for $2^4 - 1 = 15$ candidate combinations of moments. We expect to see that each combination, including (5.2), should have a posterior close to zero for large $N$ and combination $[(5.3), (5.4), (5.5)]$ should have the highest posterior probability.

TABLE 3
*Posterior probabilities, $\sigma_n^2 = n$. Set of true moments $= \{(5.3), (5.4), (5.5)\}$*

| Moments | (5.3), (5.4), (5.5) | (5.3), (5.4) | (5.3), (5.5) | (5.4), (5.5) | (5.3) | (5.4) | (5.5) |
|---|---|---|---|---|---|---|---|
| $n=100$ | 0.0076 | 0.0863 | 0.0398 | 0.0222 | 0.3564 | 0.3304 | 0.1572 |
| $n=1000$ | 0.0546 | 0.1979 | 0.0893 | 0.0502 | 0.2568 | 0.2387 | 0.1125 |
| $n=5000$ | 0.1645 | 0.2711 | 0.1230 | 0.0677 | 0.1580 | 0.1466 | 0.0691 |

The simulation result shows that if $\sigma_n^2 = n$, then although the posterior of combinations including (5.2) goes to zero quickly, the posterior probability of $[(5.3), (5.4), (5.5)]$ is still quite small, even with $N = 5000$, and is not the largest one among other true combinations (Table 3). Hence, the choice $\sigma_n^2 = n$ is too conservative. However, when $\sigma_n^2 = n^2$, the simulation result is exactly as expected. For $N = 1000$ and 5000, the combination $[(5.3), (5.4), (5.5)]$ has the largest posterior probability (Table 4). The combinations not listed are those including (5.2). They all have almost zero posterior, as desired.

**6. Discussion.** In this paper, we assume that the interior of the identified region $\text{int}(\Omega)$ is not empty. The case when $\text{int}(\Omega)$ is empty is more complicated since there is no open set contained by $\Omega$. When $\Omega$ has no interior, moment inequality models may contain exact moment conditions:

$$Em_{1j}(X, \theta_0) \geq 0, \qquad j = 1, \ldots, r,$$
$$Em_{2j}(X, \theta_0) = 0, \qquad j = 1, \ldots, p.$$

The identified region is then defined by

$$\Omega = \{\theta : Em_1(X, \theta) \geq 0, Em_2(X, \theta) = 0\}.$$

TABLE 4
*Posterior probabilities, $\sigma_n^2 = n^2$. Set of true moments $= \{(5.3), (5.4), (5.5)\}$*

| Moments | (5.3), (5.4), (5.5) | (5.3), (5.4) | (5.3), (5.5) | (5.4), (5.5) | (5.3) | (5.4) | (5.5) |
|---|---|---|---|---|---|---|---|
| $n=100$ | 0.2344 | 0.2879 | 0.1290 | 0.0682 | 0.1192 | 0.1104 | 0.0509 |
| $n=1000$ | 0.8286 | 0.0952 | 0.0428 | 0.0241 | 0.0039 | 0.0036 | 0.0017 |
| $n=5000$ | 0.9615 | 0.0223 | 0.0101 | 0.0056 | 0.0002 | 0.0002 | 0.0001 |

One of the problems one needs to take into account when considering the asymptotic behaviors of the posterior distribution is that $\Omega$ has zero Lebesgue measure, due to the loss of dimensionality. Thus, integrating over $\Omega$ always produces zero. For reasons of brevity, we do not provide a detailed discussion of this case. We point out that in this case, a dense subset in $\Omega$ still plays an important role in characterizing the large-sample behaviors of the posterior distribution. Define $\Xi = \{\theta \in \Omega : Em_1(X,\theta) > 0\}$. By assuming that $\Xi$ is dense in $\Omega$, it can still be shown that there is a large "gap" between the large-sample posterior behavior inside and outside the identified region. Inside $\Xi$, instead of being bounded below by a positive constant, we can show that the posterior density function is bounded by a polynomial rate. However, it still goes to zero exponentially fast outside the identified region. Interested readers are referred to our technical report Liao and Jiang (2008).

In partially identified models, there are two different ways to make inferences: one is studying the identified region (including consistent estimation and constructing confidence regions), while the other is directly studying the true parameter. The simulation results demonstrate that when dealing with the first goal, a flat prior is appropriate; to achieve the second goal, an informative prior is preferable. Hence, in this case, one should include as much information on the prior as possible. We believe our Bayesian method is more advantageous than the frequentist method when dealing with the second goal since the posterior distribution can provide more information about the inside of the identified region because of the prior distribution. The simulation results have verified our beliefs.

Recently, Moon and Schorfheide (2009) have considered the Bayesian approach to partially identified models when the model can involve three types of parameters: the structural parameters of interest $\theta$, a reduced-form parameter vector $\phi$ that is point-identified by data and also a vector of auxiliary parameters $\alpha$ which links structural and reduced-form parameters via some known function $\theta = \theta(\phi, \alpha)$. For a particular value of $\phi$, the auxiliary parameter takes its value in some set $\mathcal{A}_\phi$ and the identified set can then be written as

$$\Theta(\phi) = \{\theta = \theta(\phi, \alpha) : \alpha \in \mathcal{A}_\phi\}.$$

After specifying a prior distribution for both $\phi$ and $\alpha$, and combining with a likelihood function of $\phi$, a joint posterior of $\alpha$ and $\phi$ is derived, which also determines the posterior of $\theta$ via $\theta = \theta(\phi, \alpha)$. The authors also derive the Bayesian credible sets and compare them with frequentist confidence intervals for a number of particular models where $\theta(\phi, \alpha)$ is linear in $(\phi, \alpha)$ and does not involve other functions of the unknown data distribution. However, one of the main challenges of their approach is that it often requires reparametrizations between $(\theta, \alpha)$ and $(\phi, \alpha)$. Initially, it is often more natural to place a prior on the structural parameter $\theta$ and $\alpha|\theta$, but it may be inconvenient to derive the distribution $p(\phi)$ and $p(\alpha|\phi)$ from $p(\theta, \alpha)$. Another challenge is that in some models that define the relation $\theta(\phi, \alpha)$ implicitly, if $\dim(\phi) > \dim(\theta)$, it is nontrivial to specify a prior distribution $p(\phi)$ and $p(\alpha|\phi)$ such that there is a solution $\theta = \theta(\phi, \alpha)$. Also, if $\theta(\phi, \alpha)$ involves an unknown distribution of the data-generating process, there is extra variance to account for when estimating it.

In contrast to Moon and Schorfheide (2009), we proceed in a different framework of moment inequalities, one which does not require modelling the likelihood function. We construct the posterior distribution of the structural parameter using the limited information likelihood and then study the frequentist properties of the posterior. In addition, we also study the problem of model/moment selection, which is not addressed by Moon and Schorfheide (2009).

Based on the posterior distribution, we can, in principle, construct a credible set for the true parameter conditional on the data with a required coverage probability using our method (this is beyond the scope of this paper, but it is straightforward, using the posterior density function). Moon and Schorfheide (2009) have derived a Bayesian credible set for the true parameter and then compared it with the frequentist confidence interval and concluded that while frequentist confidence intervals usually extend beyond the boundaries of the identified set, the Bayesian credible sets tend to be located in the interior of the identified set. In the framework of this paper, it is also possible to derive a Bayesian credible set for the identified region if one can express the identified region explicitly in terms of $\theta$ and $\lambda$, an interesting topic for future work.

## APPENDIX A: PROOFS FOR SECTION 2

**A.1. Proof of Theorem 2.1.** Let $g(\Omega)^{-\varepsilon} = \{x \in g(\Omega) : d(x, g(\Omega)^c) \geq \varepsilon\}$, $g(\Omega)^{+\varepsilon} = \{x \in g(\Theta) : d(x, g(\Omega)) \leq \varepsilon\}$.

For all $\varepsilon > 0$, we proceed in two steps: first, show $\exists N \in \mathbb{N}$ such that when $n > N$, $\forall \varepsilon > 0$,

$$g(\Omega)^{-\varepsilon} \subset \hat{g}$$

and then show $\exists N \in \mathbb{N}$ such that when $n > N$, $\forall \varepsilon > 0$, $\hat{g} \subset g(\Omega)^{+\varepsilon}$.

Let $\inf g(\Omega) = \inf_{\theta \in \Omega} g(\theta)$ and $\sup g(\Omega) = \sup_{\theta \in \Omega} g(\theta)$.

Step I-1. Show that $g(\Omega) = [\inf g(\Omega), \sup g(\Omega)]$: obviously, $g(\Omega) \subset [\inf g(\Omega), \sup g(\Omega)]$. On the other hand, $\forall x \in [\inf g(\Omega), \sup g(\Omega)]$, since $\Omega$ is compact, $\exists \theta_1, \theta_2 \in \Omega$ so that $g(\theta_1) \leq x \leq g(\theta_2)$. By assumptions, $\Omega$ is connected and $g$ is continuous. By the intermediate value theorem, $\exists \theta^* \in \Omega$, $x = g(\theta^*)$. Hence, $x \in g(\Omega)$.

Step I-2. Show that $\exists \theta^* \in A$ and a ball $B(\theta^*, R^*)$ such that $B(\theta^*, R^*) \subset \{\theta \in \Theta : g(\theta) \leq \inf_{\theta \in \Omega} g(\Omega)^{-\varepsilon}\}$: in fact, $\forall \varepsilon > 0$, it follows by step I-1 that $g(\Omega)^{-\varepsilon} = [\inf g(\Omega) + \varepsilon, \sup g(\Omega) - \varepsilon]$. Hence, $\inf_{\theta \in \Omega} g(\Omega)^{-\varepsilon} = \inf g(\Omega) + \varepsilon$. Moreover, $\exists \theta_1 \in \Omega$, $g(\theta_1) < \inf g(\Omega) + \varepsilon$. By the continuity of $g$, there exists a ball $B(\theta_1, R)$ such that $\forall \omega \in B(\theta_1, R)$, $g(\omega) < \inf g(\Omega) + \varepsilon$. Hence, $B(\theta_1, R) \subset \{\theta \in \Theta : g(\theta) \leq \inf_{\theta \in \Omega} g(\Omega)^{-\varepsilon}\}$.

If $\theta_1 \in A$, then let $\theta^* = \theta_1$, $R^* = R$. If $\theta_1 \in \Omega \setminus A$, since $A$ is dense in $\Omega$, $B(\theta_1, \frac{R}{2}) \cap A \neq \phi$. Arbitrarily pick up an element $\theta_2 \in A \cap B(\theta_1, \frac{R}{2})$, $\forall \theta \in B(\theta_2, \frac{R}{4})$, then $d(\theta, \theta_1) \leq d(\theta, \theta_2) + d(\theta_2, \theta_1) \leq \frac{R}{4} + \frac{R}{2} < R$. Hence, $\theta \in B(\theta_1, R)$. It follows that $B(\theta_2, \frac{R}{4}) \subset B(\theta_1, R) \subset \{\theta \in \Theta : g(\theta) \leq \inf_{\theta \in \Omega} g(\Omega)^{-\varepsilon}\}$ and $\theta_2 \in A$. Let $\theta^* = \theta_2$, $R^* = \frac{R}{4}$.

Step I-3. Show that $g(\Omega)^{-\varepsilon} \subset \hat{g}$ for large $n$: by assumption 2 of Theorem 2.1, for $\theta^*$, there exists $R_{\theta^*}$ and $N \in \mathbb{N}$ such that when $\rho < R_{\theta^*}$ and $n > N$, $P(\theta \in B(\theta^*, \rho)|X^n) > \pi_n$ w.p.a.1. If we let $R_1 = \min\{R_{\theta^*}, R^*\}$, then $B(\theta^*, R_1) \subset \{\theta \in \Theta : g(\theta) \leq \inf_{\theta \in \Omega} g(\Omega)^{-\varepsilon}\}$. Hence, when $n > N$, $\forall x \in g(\Omega)^{-\varepsilon}$,

$$F_g(x) = P(g(\theta) \leq x|X^n) \geq P(g(\theta) \leq \inf g(\Omega)^{-\varepsilon}|X^n)$$
$$\geq P(\theta \in B(\theta^*, R_1)|X^n) > \pi_n.$$

Hence, $x \geq F_g^{-1}(\pi_n)$. Likewise, we can show that $x \leq F^{-1}(1 - \pi_n)$. Therefore, $g(\Omega)^{-\varepsilon} \subset [F_g^{-1}(\pi_n), F_g^{-1}(1 - \pi_n)]$.

Step II. Show for large $n$ that $\hat{g} \subset g(\Omega)^{+\varepsilon}$: step I-1 implies that $g(\Omega)^{+\varepsilon} = [\inf g(\Omega) - \varepsilon, \sup g(\Omega) + \varepsilon]$. $\forall x \in [g(\Omega)^{+\varepsilon}]^c$, either $x < \inf g(\Omega) - \varepsilon$ or $x > \sup g(\Omega) + \varepsilon$. If $x < \inf g(\Omega) - \varepsilon$, then $\{\theta \in \Theta : g(\theta) \leq x\} \subset \{\theta \in \Theta : g(\theta) \leq \inf g(\Omega) - \varepsilon\}$. In addition, since $g$ is continuous on $\Theta$, $\exists \delta > 0$ such that when $d(\theta, \Omega) \leq \delta$, $g(\theta) > \inf g(\Omega) - \varepsilon$. Therefore, $\forall \theta \in \{\theta : g(\theta) \leq \inf g(\Omega) - \varepsilon\}$, $d(\theta, \Omega) > \delta$, which implies that $\{\theta : g(\theta) \leq \inf g(\Omega) - \varepsilon\} \subset (\Omega^c)^{-\delta}$. By assumption 1 of Theorem 2.1, $\exists N \in \mathbb{N}$ such that when $n > N$, $P(\theta \in (\Omega^c)^{-\delta}|X^n) < \pi_n$ w.p.a.1. It follows that

$$P(g(\theta) \leq x|X^n) \leq P(g(\theta) \leq \inf g(\Omega) - \varepsilon|X^n) \leq P(\theta \in (\Omega^c)^{-\delta}|X^n) < \pi_n.$$

Hence, $x \leq F_g^{-1}(\pi_n)$. If $x > \sup g(\Omega) + \varepsilon$, then, by a similar argument, we can show that $x \geq F_g^{-1}(1 - \pi_n)$. Therefore, for $n > N$, if $x \in [F_g^{-1}(\pi_n), F_g^{-1}(1 - \pi_n)]$, then $x \in g(\Omega)^{+\varepsilon}$. This implies that $\hat{g} \subset g(\Omega)^{+\varepsilon}$.

Combining steps I and II, since $\varepsilon$ is arbitrary, $d_H(\hat{g}, g(\Omega)) \to 0$ in probability.

## APPENDIX B: PROOFS FOR SECTION 3

Throughout the proofs, $\phi$ denotes the empty set and $\mu(A)$ denotes the Lebesgue measure of set $A$.

**B.1. Proof of Lemma 3.1.** Recall that $(\Omega^c)^{-\varepsilon} = \{\theta : d(\theta, \Omega) \geq \varepsilon\}$, which is compact. $\forall \theta \in (\Omega^c)^{-\varepsilon}$, $\min_j \frac{Em_j(X,\theta)}{\sqrt{v_{jj}}} < 0$. $\exists \theta^* \in (\Omega^c)^{-\varepsilon}$ such that $\sup_{\theta \in (\Omega^c)^{-\varepsilon}} \min_j \frac{Em_j(X,\theta)}{\sqrt{v_{jj}}} = \min_j \frac{Em_j(X,\theta^*)}{\sqrt{v_{jj}}} < 0$. If we let

$$\delta = - \sup_{\theta \in (\Omega^c)^{-\varepsilon}} \min_j \frac{Em_j(X,\theta)}{\sqrt{v_{jj}}} > 0,$$

then $\forall \theta \in (\Omega^c)^{-\varepsilon}$, $\min_j \frac{Em_j(X,\theta)}{\sqrt{v_{jj}}} \leq -\delta < -\frac{\delta}{2}$, which implies that $(\Omega^c)^{-\varepsilon} \subset A_{\delta/2}$. Hence, $P(\theta \in (\Omega^c)^{-\varepsilon} | X^n) \leq P(\theta \in A_{\delta/2} | X^n) = o_p(a_n)$.

**B.2. Proof of Theorem 3.1.** The following lemma is useful.

LEMMA B.1. *With probability 1,*

(B.1) $\qquad P(Z \geq 0) \geq 1 - p \cdot \Phi\left(-\sqrt{n} \min_j \left\{\frac{\bar{m}_j(\theta) - (V\psi)_j/n}{\sqrt{v_{jj}}}\right\}\right),$

(B.2) $\qquad P(Z \geq 0) \leq \Phi\left(\sqrt{n} \min_j \left\{\frac{\bar{m}_j(\theta) - (V\psi)_j/n}{\sqrt{v_{jj}}}\right\}\right).$

PROOF. Let $Z = (Z_1, \ldots, Z_p)^T$.

(B.1): $P(Z \geq 0) = 1 - P\left(\bigcup_{j \leq p} Z_j < 0\right) \geq 1 - \sum_{j=1}^p P(Z_j < 0)$

$$\geq 1 - \sum_{j=1}^p \Phi\left(-\sqrt{n} \frac{\bar{m}_j(\theta) - (V\psi)_j/n}{\sqrt{v_{jj}}}\right)$$

$$\geq 1 - p \cdot \Phi\left(-\sqrt{n} \min_j \left\{\frac{\bar{m}_j(\theta) - (V\psi)_j/n}{\sqrt{v_{jj}}}\right\}\right).$$

(B.2): $P(Z \geq 0) \leq \min_j P(Z_j \geq 0) = \Phi\left(\sqrt{n} \min_j \left\{\frac{\bar{m}_j(\theta) - (V\psi)_j/n}{\sqrt{v_{jj}}}\right\}\right).$ □

PROOF OF THEOREM 3.1.1. According to Lemma 3.1, it suffices to show that, w.p.a.1, for any $\delta > 0$, $P(\theta \in A_\delta | X^n) = o_p(e^{-\alpha n})$ for some $\alpha > 0$. Define

$$\hat{A}_\delta = \left\{\theta : \min_j \frac{\bar{m}_j(X,\theta)}{\sqrt{v_{jj}}} < -\delta\right\}.$$

Then

$$P(\theta \in A_\delta | X^n) \propto \int_{A_\delta} p(\theta) L(\theta) \, d\theta$$

$$= \int_{A_\delta \cap \hat{A}_\delta} p(\theta) L(\theta) \, d\theta + \int_{A_\delta \cap \hat{A}_\delta^c} p(\theta) L(\theta) \, d\theta$$

$$\leq \int_{\hat{A}_\delta} p(\theta) L(\theta) \, d\theta + \int_{A_\delta \cap \hat{A}_\delta^c} p(\theta) L(\theta) \, d\theta,$$

$$A_\delta \cap \hat{A}_\delta^c = \left\{ \theta : \min_j \frac{Em_j(X,\theta)}{\sqrt{v_{jj}}} < -\delta \right\}$$

$$\cap \left\{ \theta : \frac{\bar{m}_i(X,\theta)}{\sqrt{v_{ii}}} \geq -\delta, i = 1, \ldots, p \right\}$$

$$= \left( \bigcup_{j=1}^{p} \left\{ \theta : \frac{Em_j(X,\theta)}{\sqrt{v_{jj}}} < -\delta \right\} \right)$$

$$\cap \left\{ \theta : \frac{\bar{m}_i(X,\theta)}{\sqrt{v_{ii}}} \geq -\delta, i = 1, \ldots, p \right\}$$

$$= \bigcup_{j=1}^{p} \left( \left\{ \theta : \frac{Em_j(X,\theta)}{\sqrt{v_{jj}}} < -\delta \right\} \right.$$

$$\left. \cap \left\{ \theta : \frac{\bar{m}_i(X,\theta)}{\sqrt{v_{ii}}} \geq -\delta, i = 1, \ldots, p \right\} \right)$$

$$= \bigcup_{j=1}^{p} A_j,$$

where

$$A_j = \left\{ \theta : \frac{Em_j(X,\theta)}{\sqrt{v_{jj}}} < -\delta \right\} \cap \left\{ \theta : \frac{\bar{m}_i(X,\theta)}{\sqrt{v_{ii}}} \geq -\delta, i = 1, \ldots, p \right\}.$$

By the weak law of large numbers, $A_j \to \phi$. Hence, $\mu(A_j) = 0$ for any $j$. Then $\mu(A_\delta \cap \hat{A}_\delta^c) = \mu(\bigcup_j A_j) \leq \sum_j \mu(A_j) = 0$ w.p.a.1. Thus, w.p.a.1, $P(\theta \in A_\delta | X^n) \leq \text{Const} \int_{\hat{A}_\delta} p(\theta) L(\theta) \, d\theta$. In addition, w.p.a.1, for some $\varepsilon > 0$,

$$L(\theta) = P(Z \geq 0) e^{-\psi^T \bar{m}(\theta) + 1/(2n) \psi^T V \psi} \prod_i \psi_i$$

$$\leq \text{Const} \cdot P(Z \geq 0) e^{\|\psi\|(\sup_{\theta \in \Theta} \|Em(X,\theta)\| + \varepsilon) + \varepsilon}$$

$$\leq \text{Const} \cdot \Phi\left(\sqrt{n}\min_j \frac{\bar{m}_j(X,\theta)}{\sqrt{v_{jj}}} + O_p\left(\frac{1}{\sqrt{n}}\right)\right).$$

Therefore, w.p.a.1,

$$P(\theta \in A_\delta | X^n) \leq \text{Const} \cdot \int_{\hat{A}_\delta} p(\theta)\Phi\left(\sqrt{n}\min_j \frac{\bar{m}_j(X,\theta)}{\sqrt{v_{jj}}} + O_p\left(\frac{1}{\sqrt{n}}\right)\right) d\theta$$

$$\leq \text{Const} \cdot \Phi\left(-\delta\sqrt{n} + O_p\left(\frac{1}{\sqrt{n}}\right)\right)$$

$$\leq \text{Const} \cdot \Phi\left(-\frac{\delta}{2}\sqrt{n}\right)$$

$$= o_p(e^{-\delta^2/8n}).$$

2. Define

$$\Omega_n = \left\{\theta : \min_j \frac{\bar{m}_j(\theta)}{\sqrt{v_{jj}}} > 0\right\}.$$

By Fatou's lemma, w.p.a.1,

$$\liminf_{n\to\infty} \int_\Xi p(\theta)L(\theta)\,d\theta$$

$$\geq \int_\Xi \liminf_{n\to\infty} p(\theta)L(\theta)\,d\theta \geq \int_{\Xi\cap\Omega\cap\Omega_n} \liminf_{n\to\infty} p(\theta)L(\theta)\,d\theta$$

$$\geq \int_{\Xi\cap\Omega\cap\Omega_n} p(\theta) \liminf_{n\to\infty} \left(1 - p\cdot\Phi\left(-\sqrt{n}\min_j \frac{\bar{m}_j(\theta) - (V\psi)_j/n}{\sqrt{v_{jj}}}\right)\right) d\theta$$

$$\geq \int_{\Xi\cap\Omega\cap\Omega_n} p(\theta)\left(1 - p\limsup_{n\to\infty}\cdot\Phi\left(-\sqrt{n}\min_j \frac{\bar{m}_j(\theta)}{\sqrt{v_{jj}}} + O_p\left(\frac{1}{\sqrt{n}}\right)\right)\right) d\theta$$

$$\geq \int_{\Xi\cap\Omega\cap\Omega_n} p(\theta)\left(1 - p\limsup_{n\to\infty}\cdot\Phi\left(-\sqrt{n}\min_j \frac{\bar{m}_j(\theta)}{2\sqrt{v_{jj}}}\right)\right) d\theta$$

$$\geq \inf_{\theta\in\Omega} p(\theta)\mu(\Xi\cap\Omega\cap\Omega_n),$$

$$\mu(\Xi\cap\Omega\cap\Omega_n)$$
$$= \mu(\Xi\cap\Omega) - \mu(\Xi\cap\Omega\cap\Omega_n^c) = \mu(\Xi) - \mu(\Xi\cap\Omega\cap\Omega_n^c)$$
$$\geq \mu(\Xi) - \mu(\Omega\cap\Omega_n^c),$$

where

$$\Omega\cap\Omega_n^c = \Omega\cap\left(\bigcup_j\left\{\theta:\frac{\bar{m}_j(\theta)}{\sqrt{v_{jj}}}\leq 0\right\}\right) = \bigcup_j\left(\Omega\cap\left\{\theta:\frac{\bar{m}_j(\theta)}{\sqrt{v_{jj}}}\leq 0\right\}\right).$$

Therefore, w.p.a.1,

$$\mu(\Omega \cap \Omega_n^c) \leq \sum_{j=1}^{p} \mu\left(\Omega \cap \left\{\theta : \frac{\bar{m}_j(\theta)}{\sqrt{v_{jj}}} \leq 0\right\}\right) \leq \sum_{j=1}^{p} \mu(\theta : Em_j(X, \theta) = 0) = 0.$$

It follows that $\mu(\Xi \cap \Omega \cap \Omega_n) \geq \mu(\Xi) > 0$. Since $p(\theta)$ is also bounded away from zero on $\Omega$, $\liminf_{n \to \infty} P(\theta \in \Xi | X^n) > 0$ in probability. $\square$

**B.3. Proof of Theorem 3.2.** In Theorem 2.1, let $A = \text{int}(\Omega)$, dense in $\Omega$. $\forall \omega \in \text{int}(\Omega), \exists R > 0$ such that $B(\omega, R) \subset \Omega$. Since $\pi_n \to 0$ but $P(\theta \in B(\omega, R) | X^n)$ is bounded away from 0 according to part 2 of Theorem 3.1, we have that for large $n$, $P(\theta \in B(\omega, R) | X^n) > \pi_n$. Therefore, by Theorem 2.1,

$$[F_g^{-1}(\pi_n), F_g^{-1}(1 - \pi_n)] \to g(\Omega) \qquad \text{in probability.}$$

**B.4. Proof of Theorem 3.3.** To show this theorem, the following lemmas are useful.

LEMMA B.2. *In probability,*

(B.3) $$\limsup_{n \to \infty} \max_{\theta \in \Theta} \ln p(\theta | X^n) < \infty.$$

$\forall \varepsilon > 0$,

(B.4) $$\liminf_{n \to \infty} \inf_{\theta \in \Omega^{-\varepsilon}} p(\theta | X^n) > 0.$$

PROOF. (B.3): For some $\varepsilon > 0$,

$$\limsup_{n \to \infty} \sup_{\theta \in \Theta} L(\theta) \leq \prod_j \psi_j e^{\|\psi\|(\sup_{\theta \in \Theta} \|Em(X,\theta)\| + \varepsilon) + \varepsilon} < \infty.$$

Thus,

$$\limsup_{n \to \infty} \max_{\theta \in \Theta} \ln p(\theta | X^n) = \text{Const} \cdot \limsup_{n \to \infty} \max_{\theta \in \Theta} \ln p(\theta) L(\theta)$$

$$\leq C \cdot \ln\left(\sup_{\theta \in \Theta} p(\theta) \cdot \limsup_{n \to \infty} \sup_{\theta \in \Theta} L(\theta)\right) < \infty.$$

(B.4): $\forall \varepsilon > 0$,

$$\liminf_{n \to \infty} \inf_{\theta \in \Omega^{-\varepsilon}} L(\theta) \geq \text{Const} \cdot \liminf_{n \to \infty} \inf_{\theta \in \Omega^{-\varepsilon}} P(Z_\theta \geq 0) e^{-\|\psi\| \cdot (\sup_{\theta \in \Theta} \|Em(X,\theta)\| + \varepsilon)}$$

$$\geq C \cdot \liminf_{n \to \infty} \inf_{\theta \in \Omega^{-\varepsilon}} P(Z_\theta \geq 0) > 0.$$

Here, $C$ denotes a positive constant. The last inequality follows since $Z_\theta \sim N_p(\bar{m}(\theta) - V\psi/n, V/n)$, $\Omega^{-\varepsilon} \subset \Omega$ and $Em(X, \theta) \geq 0$ on $\Omega$. $\square$

LEMMA B.3. *In probability,*

1. *for all $\varepsilon > 0$,*

$$\limsup_{n\to\infty} \sup_{\theta\in\Omega^{-\varepsilon}} |\max_{\omega\in\Theta} \ln p(\omega|X^n) - \ln p(\theta|X^n)| < \infty;$$

2. *if $\varepsilon_n \prec n$, then $\forall \varepsilon > 0$,*

$$\frac{\varepsilon_n}{\inf_{\theta\in(\Omega^c)^{-\varepsilon}} |\ln p(\theta|X^n)|} \to 0.$$

PROOF. 1. For each $n$,

$$\sup_{\theta\in\Omega^{-\varepsilon}} |\max_{\omega\in\Theta} \ln p(\omega|X^n) - \ln p(\theta|X^n)| = \max_{\theta\in\Theta} \ln p(\theta|X^n) - \inf_{\theta\in\Omega^{-\varepsilon}} \ln p(\theta|X^n).$$

The result follows immediate from Lemma B.2.

2. W.p.a.1, $\ln p(\theta|X^n) < 0$ on $(\Omega^c)^{-\varepsilon}$, hence

$$\inf_{\theta\in(\Omega^c)^{-\varepsilon}} |\ln p(\theta|X^n)|$$

$$= - \sup_{\theta\in(\Omega^c)^{-\varepsilon}} \ln p(\theta|X^n)$$

$$\geq - \mathrm{Const} \cdot \ln \sup_{\theta\in(\Omega^c)^{-\varepsilon}} L(\theta)$$

$$\geq -C \cdot \ln \sup_{\theta\in(\Omega^c)^{-\varepsilon}} P(Z_\theta \geq 0)$$

$$\geq -C \cdot \ln \sup_{\theta\in(\Omega^c)^{-\varepsilon}} \Phi\left(\sqrt{n}\min_j \frac{\bar{m}_j(\theta) - (V\psi)_j/n}{\sqrt{v_{jj}}}\right).$$

As shown in the proof of Lemma 3.1, there exists some $\delta > 0$ such that $(\Omega^c)^{-\varepsilon} \subset A_\delta$, where $A_\delta = \{\theta : \min_j \frac{Em_j(X,\theta)}{\sqrt{v_{jj}}} < -\delta\}$. Thus, w.p.a.1,

$$\inf_{\theta\in(\Omega^c)^{-\varepsilon}} |\ln p(\theta|X^n)| \geq -C \cdot \ln \sup_{\theta\in A_\delta} \Phi\left(\sqrt{n}\min_j \frac{\bar{m}_j(\theta) - (V\psi)_j/n}{\sqrt{v_{jj}}}\right)$$

$$\geq -C \cdot \ln \sup_{\theta\in A_\delta} \Phi\left(\sqrt{n}\min_j \frac{\bar{m}_j(\theta)}{2\sqrt{v_{jj}}}\right)$$

$$\geq -C \cdot \ln \Phi\left(-\frac{\delta}{2}\sqrt{n}\right)$$

$$\geq -C_1 \cdot n + C_2 \ln n + C_3,$$

where $C_1 > 0, C_2$ and $C_3$ denote finite constants. This implies that $\inf_{\theta\in(\Omega^c)^{-\varepsilon}} |\ln p(\theta|X^n)| = O_p(n)$.

□

PROOF OF THEOREM 3.3. For all $\varepsilon > 0$, since $\varepsilon_n \to \infty$, we have, by part 1 of Lemma B.3, that $\exists N \in \mathbb{N}$ such that when $n > N$, for any $\theta \in \Omega^{-\varepsilon}$,

$$\max_{\omega \in \Theta} \ln p(\omega|X^n) - \ln p(\theta|X^n) < \varepsilon_n, \qquad \text{w.p.a.1.}$$

Therefore, when $n > N$, $\Omega^{-\varepsilon} \subset A_n$, which implies that $\limsup_{n \to \infty} \sup_{\theta \in \Omega} d(\theta, A_n) \leq \varepsilon$.

On the other hand, let $M = \liminf_{n \to \infty} \max_{\theta \in \Theta} \ln p(\theta|X^n)$. By (B.3) in Lemma B.2, $M < \infty$. Moreover, by (B.4),

$$M \geq \liminf_{n \to \infty} \inf_{\theta \in \Omega^{-\varepsilon}} \ln p(\theta|X^n) \geq \ln \liminf_{n \to \infty} \inf_{\theta \in \Omega^{-\varepsilon}} p(\theta|X^n) > -\infty.$$

Hence, $M \in \mathbb{R}$ and, by the definition of $M$, $\exists N_1 \in \mathbb{N}$ such that when $n > N_1$,

$$\max_{\theta \in \Theta} \ln p(\theta|X^n) > M - \varepsilon.$$

In addition, $\forall \theta \in (\Omega^c)^{-\varepsilon}$, $p(\theta|X^n) \to 0$ in probability. Thus, for large $n$, $\ln p(\theta|X^n) < 0$ on $(\Omega^c)^{-\varepsilon}$. $\exists N_2 \in \mathbb{N}$ such that when $n > N_2$,

$$\inf_{\theta \in (\Omega^c)^{-\varepsilon}} |\ln p(\theta|X^n)| = -\sup_{\theta \in (\Omega^c)^{-\varepsilon}} \ln p(\theta|X^n) > \varepsilon_n - (M - \varepsilon),$$

where the inequality follows by part 2 of Lemma B.3. Therefore, when $n > N_2$,

(B.5) $$\sup_{\theta \in (\Omega^c)^{-\varepsilon}} \ln p(\theta|X^n) < -\varepsilon_n + (M - \varepsilon).$$

However, when $n > \max\{N_1, N_2\}$, $\forall \theta \in A_n = \{\theta : \max_{\omega \in \Theta} \ln p(\omega|X^n) - \ln p(\theta|X^n) \leq \varepsilon_n\}$, $\ln p(\theta|X^n) \geq \max_{\omega \in \Theta} \ln p(\omega|X^n) - \varepsilon_n > M - \varepsilon - \varepsilon_n$. Comparing this with (B.5), we see that $\theta \notin (\Omega^c)^{-\varepsilon}$. In other words, $d(\theta, \Omega) < \varepsilon$. It follows that

$$\limsup_{n \to \infty} \sup_{\theta \in A_n} d(\theta, \Omega) \leq \varepsilon.$$

Since $\varepsilon$ is arbitrary, $d_H(A_n, \Omega) \to 0$ in probability. □

## APPENDIX C: PROOFS FOR SECTION 4

LEMMA C.1. *If we suppose that $C_s = (M_s, \Theta_s)$, $\Omega_s = \{\theta \in \Theta_s : EM_s(X, \theta) \geq 0\}$ and $\Theta_s$ is compact, then, for some $\xi \in \Theta_s$ and normalization parameter $C$,*

(C.1) $$p \lim_{n \to \infty} p(C_s|X^n) = C\left(\prod_{j \in S}\right) p(C_s) P(\theta \in \Omega_s|C_s) e^{-\psi^T EM_s(X, \xi)}.$$

PROOF. By the integral intermediate value theorem, the right-hand side of (C.1) can be written as

$$\text{(C.2)} \qquad \text{RHS} = C \left( \prod_{j \in S} \right) p(C_s) \int_{\Theta_s} 1_{\Omega_s} p(\theta|C_s) e^{-\psi^T E M_s(X,\theta)} \, d\theta.$$

On the other hand,

$$\text{(C.3)} \quad \begin{aligned} p(C_s|X^n) &= C \int \int_{\Theta_s \times [0,\infty)^m} \frac{1}{\sqrt{\det((2\pi)/nV)}} \\ &\quad \times e^{-n/2(\bar{M}_s(\theta)-\lambda)^T V^{-1}(\bar{M}_s(\theta)-\lambda)} \\ &\quad \times \left( \prod_{j \in S} \psi_j \right) e^{-\psi^T \lambda} p(\theta|C_s) p(C_s) \, d\theta \, d\lambda \\ &= C \left( \prod_{j \in S} \right) p(C_s) \int_{\Theta_s} p(\theta|C_s) p(Z_\theta \geq 0) \\ &\quad \times e^{-\psi^T \bar{M}_s(\theta) + 1/(2n)\psi^T V \psi} \, d\theta, \end{aligned}$$

where $Z_\theta \sim N_m(\bar{M}_s(\theta) - \frac{V\psi}{n}, \frac{V}{n})$. Take the difference between (C.2) and (C.3):

$$\text{(C.4)} \quad \begin{aligned} &|p(C_s|X^n) - \text{RHS}| \\ &\leq \text{Const} \cdot p(C_s) \\ &\quad \times \int_{\Theta_s} p(\theta|C_s) |1_{\Omega_s} p(\theta|C_s) e^{-\psi^T E M_s(X,\theta)} \\ &\qquad - p(Z_\theta \geq 0) e^{-\psi^T \bar{M}_s(\theta) + 1/(2n)\psi^T V \psi}| \, d\theta. \end{aligned}$$

If we let $\Delta(\theta) = p(\theta|C_s) |1_{\Omega_s} p(\theta|C_s) e^{-\psi^T E M_s(X,\theta)} - p(Z_\theta \geq 0) \times e^{-\psi^T \bar{M}_s(\theta) + 1/(2n)\psi^T V \psi}|$, then (C.4) can be rewritten as

$$|p(C_s|X^n) - \text{RHS}|$$
$$\leq \text{Const} \cdot p(C_s) \left( \int_{U_1} \Delta(\theta) \, d\theta + \int_{U_2} \Delta(\theta) \, d\theta + \int_{U_3} \Delta(\theta) \, d\theta \right),$$

where

$$U_1 = \{\theta \in \Theta_s : EM_s(X,\theta) > 0\},$$
$$U_2 = \{\theta \in \Theta_s : EM_s(X,\theta) \geq 0, Em_j(X,\theta) = 0 \text{ for some } m_j \in M_s\},$$
$$U_3 = \{\theta \in \Theta_s : \text{ for some } m_j \in M_s, Em_j(X,\theta) < 0\}.$$

We next look at the integrations on $U_i$, $i=1,2,3$.

$U_1$: Note that $\Omega_s = \{\theta \in \Theta_s : EM_s(X,\theta) \geq 0\}$ and $Z_\theta \sim N_m(\bar{M}_s(\theta) - \frac{V\psi}{n}, \frac{N}{n})$. For any $\varepsilon > 0$, by the uniform weak law of large numbers, w.p.a.1, $\sup_{\theta \in U_1} |P(Z_\theta \geq 0) - 1_{\Omega_s}| < \varepsilon$. Hence, for large $n$, w.p.a.1,

$$\sup_{\theta \in U_1} |1_{\Omega_s} p(\theta|C_s) e^{-\psi^T EM_s(X,\theta)} - p(Z_\theta \geq 0) e^{-\psi^T \bar{M}_s(\theta) + 1/(2n)\psi^T V \psi}| < \varepsilon.$$

Hence,

$$\int_{U_1} \Delta(\theta) \, d\theta \leq \varepsilon \int_{U_1} p(\theta|C_s) \, d\theta \leq \varepsilon.$$

$U_2$: The Lebesgue measure of $U_2 = 0$.

$U_3$: $\forall \theta \in U_3$, $1_{\theta \in \Omega_s} = 0$, hence,

$$\Delta(\theta) = p(\theta|C_s) P(Z_\theta \geq 0)$$
$$\times e^{-\psi^T \bar{M}_s(\theta) + 1/(2n)\psi^T V \psi} \qquad \forall \varepsilon > 0,$$

w.p.a.1, $P(Z_\theta \geq 0) < \varepsilon$, thus, for large $n$, w.p.a.1,

$$\int_{U_3} \Delta(\theta) \, d\theta \leq \varepsilon \int_{U_3} p(\theta|C_s) e^{-\psi^T \bar{M}_s(\theta)} \, d\theta$$
$$\leq e^{\|\psi\| \cdot (\sup_{\Theta_s} \|EM_s(X,\theta)\| + \varepsilon)} \varepsilon.$$

We have thus shown that $|p(C_s|X^n) - \text{RHS}| < \text{Const} \cdot p(C_s) \varepsilon$, w.p.a.1, with arbitrarily small $\varepsilon$. □

### C.1. Proof of Theorem 4.1.

1. The result follows immediately from Lemma C.1 and Assumption 4.2.
2. For some normalization parameter $C$,

$$p(C_s|X^n) = C \iint_{\Theta_s \times [0,\infty)^m} \frac{1}{\sqrt{\det(2\pi/nV)}}$$
$$\times e^{-n/2(\bar{M}_s(\theta) - \lambda)^T V^{-1}(\bar{M}_s(\theta) - \lambda)}$$
$$\times \left(\prod_{j \in S} \psi_j\right) e^{-\psi^T \lambda} p(\theta|C_s) p(C_s) \, d\theta \, d\lambda.$$

Since $V^{-1}$ is positive definite and $C_s = (M_s, \Theta_s)$ is not true, $\exists \tau > 0$ such that

$$\inf_{\Theta_s \times [0,\infty)^m} (EM_s(X,\theta) - \lambda)^T V^{-1} (EM_s(X,\theta) - \lambda) > \tau.$$

Hence, w.p.a.1, $\inf_{\Theta_s \times [0,\infty)^m}(\bar{M}_s(\theta) - \lambda)^T V^{-1}(\bar{M}_s(\theta) - \lambda) > \tau$. Therefore, w.p.a.1,

$$p(C_s|X^n) \leq C\left(\prod_{j \in S} \psi_j\right) e^{-n/2\tau}$$

$$\times \iint_{\Theta_s \times [0,\infty)^m} \frac{1}{\sqrt{\det(2\pi/nV)}}$$

$$\times e^{-\psi^T \lambda} p(\theta|C_s) p(C_s) \, d\theta \, d\lambda$$

$$\leq \text{Const} \cdot n^m e^{-n/2\tau} p(C_s) \left(\prod_S \psi_j\right) \int_{[0,\infty)^m} e^{-\psi^T \lambda} \, d\lambda$$

$$\leq e^{-\tau/4n} p(C_s).$$

**C.2. Proof of Theorem 4.2.** If we let $\Delta$ denote the index set corresponding to the moment inequalities that are selected by $C_s^2$ but *not* by $C_s^1$, then $p(C_s^2 X^n)$ has a $\prod_{j \in \Delta} \psi_j$ term that does not show up in $p(C_s^1|X^n)$. If $p(C_s^1) = p(C_s^2)$, by Lemma C.1,

$$p\lim_{n \to \infty} \text{BF}_{12} = \frac{P(\theta \in \Omega_s^1|C_s^1)}{P(\theta \in \Omega_s^2|C_s^2)} \cdot \frac{e^{-\psi_1^T E M_s^1(\xi_1)}}{e^{-\psi_2^T E M_s^2(\xi_2)}} \cdot \prod_{j \in \Delta} \frac{1}{\psi_j}.$$

Note that since $M_s^i$, $i = 1, 2$, are both true, the integral intermediate value theorem guarantees that $\xi_i \in \Omega_s^i$, hence $EM_s^i(\xi_i) \geq 0$ for $i = 1, 2$. It follows that

$$p\lim_{n \to \infty} \text{BF}_{12} \geq \frac{P(\theta \in \Omega_s^1|C_s^1)}{P(\theta \in \Omega_s^2|C_s^2)} \cdot e^{-\psi_1^T E M_s^1(\xi_1)} \cdot \prod_{j \in \Delta} \frac{1}{\psi_j}$$

$$\geq \frac{P(\theta \in \Omega_s^1|C_s^1)}{P(\theta \in \Omega_s^2|C_s^2)} \cdot e^{-\|\psi_1\| \cdot \sup_{\theta \in \Theta_s} \|Em(X,\theta)\|} \cdot \prod_{j \in \Delta} \frac{1}{\psi_j}$$

$$\geq \frac{P(\theta \in \Omega_s^1|C_s^1)}{P(\theta \in \Omega_s^2|C_s^2)}$$

$$\geq 1.$$

The third inequality is due to Assumption 4.6 and the last inequality follows from $\Omega_s^2 \subset \Omega_s^1$.

**C.3. Proof of Theorem 4.3.**

APPROACH 1. Suppose that $p(C_s) \propto n^{\alpha[\dim(M_s) - \dim(\Theta_s)]}$ for some $\alpha > 0$.

1. If $C^*$ is false, then by part 2 of Theorem 4.1, $\exists \beta > 0$ such that
$$p(C^*|X^n) = o_p(e^{-\beta n}) n^{\alpha[\dim(M_s) - \dim(\Theta_s)]},$$
which is exponentially small. However, there exists at least one true combination $C_s$, with posterior distribution bounded away from zero. Hence, w.p.a.1, $p(C_s|X^n) > p(C^*|X^n)$, a contradiction.

2. Since $C^*$ is true, $\Omega^*$, the identified region which is defined by it satisfies $p(\theta \in \Omega^*|C^*) > 0$. By Lemma C.1, $p(C^*|X^n) = O_p(n^{\alpha[\dim(M_s) - \dim(\Theta_s)]})$. It follows immediately from the definition of $C^*$ that $C^*$ has the largest value of $\dim(M_s) - \dim(\Theta_s)$.

APPROACH 2. Suppose that $p(C_s)$ is the uniform prior of $C_s$, and we put multivariate normal priors on unrestricted parameters. For any candidate $C_s$,

$$p(C_s|X^n) \propto \iiint_{\Theta_s \times [0,\infty)^m \times \mathbb{R}^{p-m}} L(X^n|\theta_s, \lambda, C_s) p(\theta_s|C_s) p(\lambda_s|C_s)$$
$$\times p(\lambda_s^c|C_s) \, d\theta_s \, d\lambda_s \, d\lambda_s^c.$$

Let

(C.5) $$L(X^n|\theta_s, \lambda_s, C_s) = \int_{\mathbb{R}^{p-m}} L(X^n|\theta_s, \lambda, C_s) p(\lambda_s^c|C_s) \, d\lambda_s^c,$$

(C.6) $$L(X^n|\theta_s, C_s) = \int_{[0,\infty)^m} L(X^n|\theta_s, \lambda_s, C_s) p(\lambda_s|C_s) \, d\lambda_s.$$

A tedious calculation shows that
$$L(X^n|\theta_s, \lambda_s, C_s) = \frac{1}{\sqrt{\det(2\pi S_n)}}$$
$$\times \exp\left\{-\frac{1}{2}(\bar{M}_s(\theta) - \lambda_s, \bar{M}_s^c(\theta)) S_n^{-1} \begin{pmatrix} \bar{M}_s(\theta) - \lambda_s \\ \bar{M}_s^c(\theta) \end{pmatrix}\right\},$$

where
$$S_n = \frac{V}{n} + \begin{pmatrix} 0 & 0 \\ 0 & \sigma_n^2 I_{p-m} \end{pmatrix}; \quad \text{we write } S_n^{-1} = n \begin{pmatrix} \Sigma_1 & \Sigma_3 \\ \Sigma_3^T & \Sigma_2 \end{pmatrix}.$$

We can then calculate (C.6):
$$L(X^n|\theta_s, C_s) = \text{Const} \cdot \frac{1}{\sqrt{\det(V_2)}} P(Z_\theta \geq 0) e^{\tau(\theta)},$$

where:

- $V_2 = V_{22} + n\sigma_n^2 I_{p-m}$, with $V_{22}$ being the lower diagonal block of $V$;
- $Z_\theta \sim N_m(\bar{M}_s(\theta) + \Sigma_1^{-1} \Sigma_3^T \bar{M}_s^c(\theta) - \frac{1}{n} \Sigma_1^{-1} \psi, \frac{\Sigma_1^{-1}}{n})$;

- $\tau(\theta) = -\frac{n}{2}\bar{M}_s^c(\theta)(V_{22} + n\sigma_n^2 I_{p-m})^{-1}\bar{M}_s^c(\theta) - \psi^T[\Sigma_1^{-1}\Sigma_3^T\bar{M}_s^c(\theta) + \bar{M}_s(\theta)] + \frac{1}{2n}\psi^T\Sigma_1^{-1}\psi$.

Given that $\|V\| = O(1)$, one can show that $\|\Sigma_1\| = O(1), \|\Sigma_2\| = O(\frac{1}{n\sigma_n^2})$ and $\|\Sigma_3\| = O(\frac{1}{n\sigma_n^2})$.

Define an operator of $n^{-1}$ and $\bar{M}(\cdot)$,

$$g(n^{-1}, \bar{M}(\theta_s)) = (n\sigma_n^2)^{(\dim(\theta_s)+p-m)/2} L(X^n|\theta_s, C_s) p(\theta_s|C_s),$$

where $L(X^n|\theta_s, C_s)$ is the integrated limited information likelihood of $(\theta_s, C_s)$, by integrating out $\lambda$. We use a factor $(n\sigma_n^2)^{(\dim(\theta_s)+p-m)/2}$ for rescaling so that

$$g(n^{-1}, \bar{M}(\theta_s)) = O_p(1).$$

Hence, w.p.a.1, without changing the orders, we have

$$g(n^{-1}, \bar{M}(\theta_s)) = (n\sigma_n^2)^{(\dim(\theta_s)+p-m)/2} L(X^n|\theta_s, C_s) p(\theta_s|C_s)$$
$$= \text{Const} \cdot P(Z_\theta \geq 0) e^{\tau(\theta)} e^{-\theta^T\theta/(2n\sigma_n^2)}.$$

This yields that

$$p(C_s|X^n) = \text{Const} \cdot \int_{\Theta_s} L(X^n|\theta_s, C_s) p(\theta_s|C_s) \, d\theta_s$$
$$= (n\sigma_n^2)^{-(\dim(\theta_s)+p-m)/2} \text{Const} \cdot \int_{\Theta_s} g(n^{-1}, \bar{M}(\theta_s)) \, d\theta_s.$$

The following lemma is needed before proceeding.

LEMMA C.2. *Under Assumption 4.6, in probability,*

$$\lim_{n\to\infty} \int_{\Theta_{s_2}} g(n^{-1}, \bar{M}(\theta_s)) \, d\theta_s = \int_{\Theta_{s_2}} \lim_{n\to\infty} g(n^{-1}, \bar{M}(\theta_s)) \, d\theta_s.$$

The proof is given at the end of this section.

Hence, by Lemma C.2, if $C_s$ is true, then, in probability,

$$\lim_{n\to\infty} \int_{\Theta_s} g(n^{-1}, \bar{M}(\theta_s)) \, d\theta_s = \int_{\Theta_s} \lim_{n\to\infty} g(n^{-1}, \bar{M}(\theta_s)) \, d\theta_s$$
$$= \text{Const} \cdot \int_{\{\theta : EM_s(\theta) \geq 0\}} e^{-\psi^T EM_s(\theta)} \, d\theta$$
$$= O(1).$$

The second equality is due to $P(Z_\theta \geq 0) \xrightarrow{p} 1_{\{\theta \in \Theta_s : EM_s(\theta) \geq 0\}}$ and $\tau(\theta) \xrightarrow{p} \psi^T EM_s(\theta)$. The last equality follows since $\{\theta \in \Theta_s : EM_s(\theta) \geq 0\}$ is the identified region of $C_s$ and is assumed to be compact.

Hence, $p(C_s|X^n) = O_p((n\sigma_n^2)^{-(\dim(\theta_s)+p-m)/2})$, which follows from the fact that the optimal $C^*$ that maximizes $p(C_s|X^n)$ has the largest value of $m - \dim(\theta_s)$.

PROOF OF LEMMA C.2. We apply the following theorem.

THEOREM C.1 [Billingsley (1986), Theorem 16.8]. *Let $f(t,w): (T \times W) \to \mathbb{R}$ be a real-valued function, absolutely integrable with respect to $w$. Suppose that:*

1. *$f(t,w)$ is continuous on a neighborhood of $t = t_0$ for almost all $w \in W$;*
2. *there exists a function $g: W \to \mathbb{R}^+$ such that $|f(t,w)| \leq g(w)$ for any $t \in T$ and $\int_W g(w)\,dw < \infty$.*

*Then*
$$\lim_{t \to t_0} \int_W f(t,w)\,dw = \int_W \lim_{t \to t_0} f(t,w)\,dw.$$

PROOF. $g(n^{-1}, \bar{M}(\theta_s)) = \text{Const} \cdot P(Z_\theta \geq 0) e^{\tau(\theta)} e^{-\theta^T \theta/(2\sigma_n^2)}$. Here, the sample moments and $\theta_s$ are separated. If we let $\bar{W}$ denote a vector of all the sample moments $\bar{W} = (\overline{ZY_1}, \overline{ZY_2}, \overline{ZX^T})$, then we can write $g(n^{-1}, \bar{M}(\theta_s)) = g(n^{-1}, \bar{W}, \theta_s)$. It suffices to show that

$$\lim_{n \to \infty} \int_{\Theta_s} g(n^{-1}, \bar{W}, \theta_s)\,d\theta_s = \int_{\Theta_s} \lim_{n \to \infty} g(n^{-1}, \bar{W}, \theta_s)\,d\theta_s.$$

We proceed by verifying the conditions in Theorem C.1.

CONDITION 1. Note that $P(Z_\theta \geq 0) \to^p 1_{\Omega_s}$ for almost all $\theta_s$, except on a zero-measure set $\{\theta: \exists j, Em_{sj}(\theta) = 0\}$. Hence, it is straightforward to verify that $g(n^{-1}, \bar{W}, \theta_s)$ is continuous on a small neighborhood of $(0, EW)$ for almost all $\theta_s$.

CONDITION 2. In this case, the moment inequality functions are all linear in $\theta$, hence we can write $\bar{M}_s^c(\theta) = a_1 + B_1 \theta$ and $\bar{M}_s(\theta) = a_2 + B_2\theta$, where $B_i$ are matrices. We can then show that, w.p.a.1 (we omit some intermediate calculations),
$$e^{\tau(\theta)} \leq Ce^{-\alpha^T/(2n\sigma_n^2)\theta_s - \beta^T \theta_s},$$
where $C > 0$ and $\alpha$, $\beta$ are constant vectors. Hence, for large $n$,

$$e^{\tau(\theta)} e^{-\theta^T \theta/(2n\sigma_n^2)} \leq Ce^{-\beta^T \theta} e^{-1/(2n\sigma_n^2)(\theta+\alpha/2)^T(\theta+\alpha/2)} e^{\alpha^T \alpha/(8n\sigma_n^2)}$$
(C.7)
$$\leq \text{Const} \cdot e^{-\beta^T \theta}.$$

Furthermore, since $\text{rank}(EZX^T) = \dim(X)$, there exists $c > 0$ such that $\forall \theta_s \in \Theta_s$, we can always find a component $Em_{sk}(\theta)$ of $EM_s(\theta)$ such that $Em_{sk}(\theta) < -c\|\theta\|$. Write $\theta_s = \omega r$, where $\omega$ and $r$ denote the unit direction vector and the radius of $\theta_s$, respectively. Then $Em_{sk}(\theta) < -cr$ (here, $k$ and $r$ depend on $\theta$, but $c$ does not). For $Z_\theta \sim N_m(EM_s(\theta) + O_p(\frac{1}{\sqrt{n}}), \frac{\Sigma_1^{-1}}{n})$, $\forall \varepsilon > 0$, w.p.a.1, for some $v_k > 0$,

$$
\begin{aligned}
P(Z_\theta \geq 0) &= P(N_m(\bar{M}_s(\theta) + O_p(n^{-1}), \Sigma_1^{-1}/n) \geq 0) \\
&\leq P(N(\bar{m}_{sk}(\theta) + O_p(n^{-1}), v_k/n) \geq 0) \\
&\leq P(N(Em_{sk}(\theta) + \varepsilon, v_k/n) \geq 0) \\
&= 1 - \Phi\left(-\sqrt{n}\frac{Em_{sk}(\theta) + \varepsilon}{\sqrt{v_k}}\right) \\
&\leq 1 - \Phi\left(\frac{\sqrt{n}(cr - \varepsilon)}{\sqrt{v_k}}\right) \\
&\leq \sqrt{\frac{v_k}{n}} \frac{1}{\sqrt{2\pi}(cr - \varepsilon)} e^{-(1/2)(n/v_k)(cr - \varepsilon)^2}.
\end{aligned}
$$
(C.8)

The last inequality follows from Mill's ratio inequality. We can choose $\varepsilon = \frac{cr}{2}$, then, for large $n$, $P(Z_\theta \geq 0) \leq \text{Const} \cdot e^{-c^2 r^2/(8v_k)}$. Combining with (C.7), we obtain an integrable function to upper bound $g(n^{-1}, \bar{M}(\theta_s))$: for all $n$,

$$
\begin{aligned}
g(n^{-1}, \bar{M}(\theta_s)) &\leq \text{Const} \cdot e^{-\beta^T \theta} e^{-c^2 r^2/(8v_k)} \\
&= \text{Const} \cdot e^{-c^2/(8v_k)(r + 4v_k \beta^T \omega/c^2)^2} e^{2v_k(\beta^T \omega)^2/c^2} \\
&\leq \text{Const} \cdot e^{-c^2/(8v_k)(r + 4v_k \beta^T \omega/c^2)^2}.
\end{aligned}
$$

To see that this upper bound function is integrable, write $\theta = r\omega$, so

$$
\int_{\Theta_s} e^{-c^2/(8v_k)(r + 4v_k \beta^T \omega/c^2)^2} d\theta_s = \int_0^\infty e^{-c^2/(8v_k)(r + 4v_k \beta^T \omega/c^2)^2} dr \oint_{\{\omega : \|\omega\| = 1\}} d\omega
$$
$$
\leq \sqrt{\frac{8\pi v_k}{c^2}} S(\{\omega : \|\omega\| = 1\}),
$$

where $S(\{\omega : \|\omega\| = 1\})$ denotes the surface area of the unit ball $\{\omega : \|\omega\| = 1\}$. □

*Completion of the Proof of Theorem 4.3, Approach 2.* It is left to show that the posterior of a false combination is exponentially small. Let $C_s = (EM_s, \Theta_s)$ be a false combination. Then, by definition, we can write $\xi = \inf_{\Theta_s \times [0, \infty)^{\dim(\lambda)}} \|EM_s(X, \theta) - \lambda\|^2 > 0$. Define a compact ball $B(d_n) = \{\theta : \|\theta\| \leq$

42                        Y. LIAO AND W. JIANG

$d_n\}$ and $U_n = \{\theta : \|\theta\| > d_n\}$ for some radius $d_n \to \infty$, with the rate to be specified below. Then

$$p(C_s|X^n) \propto \int_{B(d_n)} L(X^n|\theta, C_s) p(\theta_s|C_s)\, d\theta_s + \int_{U_n} L(X^n|\theta, C_s) p(\theta_s|C_s)\, d\theta_s.$$

On one hand, w.p.a.1, we can show that

$$L(X^n|\theta, C_s) \leq \text{Const} \cdot e^{-\xi n/2} \det(S_n)^{-1/2} e^{-n\bar{M}_s^c(\theta)^T \Sigma_3^T \bar{M}_s(\theta) - n/2 \bar{M}_s^c(\theta)^T \Sigma_2 \bar{M}_s^c(\theta)}$$

$$\times \int_{\lambda_s \geq 0} e^{-(\psi - n\Sigma_3 \bar{M}_s^c(\theta))^T \lambda_s}\, d\lambda_s$$

$$\leq \text{Const} \cdot e^{-\xi n/2} \det(S_n)^{-1/2} e^{-n\bar{M}_s^c(\theta)^T \Sigma_3^T \bar{M}_s(\theta)}.$$

Hence, $\int_{B(d_n)} L(X^n|\theta, C_s) p(\theta_s|C_s)\, d\theta_s \leq \text{Const}(n\sigma_n^2)^{-(\dim(\Theta_s) + p - m)/2} e^{-\xi n/2}$. $\int_{B(\omega, d_n)} e^{-n\bar{M}_s^c(\theta)^T \Sigma_3^T \bar{M}_s(\theta)}\, d\theta$. Note that $\|\Sigma_3\| = O(\frac{1}{n\sigma_n^2})$ and, in this example, $\bar{M}(\theta)$ is linear in $\theta$, hence $\sup_{\theta \in B(d_n)} \|n\bar{M}_s^c(\theta)^T \Sigma_3^T \bar{M}_s(\theta)\| \leq \text{Const}(d_n/\sigma_n)^2$, w.p.a.1. Assuming that $(d_n/\sigma_n)^2 = o_p(n)$, we have $\xi n \succ (\frac{d_n}{\sigma_n})^2$. Hence,

$$\int_{B(d_n)} L(X^n|\theta, C_s) p(\theta_s|C_s)\, d\theta_s$$
$$\leq (n\sigma_n^2)^{-(\dim(\Theta_s) + p - m)/2} e^{-\xi n/2} e^{c(d_n/\sigma_n)^2} d_n^{\dim(\Theta_s)}$$
$$\leq e^{-\alpha n}$$

for some $\alpha > 0$. On the other hand,

$$\int_{U_n} L(X^n|\theta, C_s) p(\theta_s|C_s)\, d\theta_s$$
$$\propto (n\sigma_n^2)^{-(\dim(\Theta_s) + p - m)/2} \int_{U_n} P(Z_\theta \geq 0) e^{\tau(\theta)} e^{-\theta^T \theta/(2n\sigma_n^2)}\, d\theta.$$

We use (C.7), $e^{\tau(\theta)} e^{-\theta^T \theta/(2n\sigma_n^2)} \leq \text{Const} \cdot e^{-\beta^T \theta}$, for some constant vector $\beta$. Combining with (C.8) and using the same trick as before by writing $\theta = \omega\|\theta\|$, we have

$$\int_{U_n} L(X^n|\theta, C_s) p(\theta_s|C_s)\, d\theta_s \leq \text{Const} \cdot \int_{d_n}^{\infty} e^{-an(r+b)^2}\, dr \cdot d_n^{\dim(\Theta_s)},$$

where $a > 0, b \in \mathbb{R}$ are constant. By Mill's ratio inequality, it is less than $e^{-\alpha n d_n^2}$ for some $\alpha > 0$. $\square$

seg**Acknowledgments.** The authors thank Ivan Canay, Joel Horowitz and Elie Tamer for many helpful comments and discussions. We also thank the Editor and an anonymous referee for useful comments that have improved the presentation of this paper.

Department of Statistics
Northwestern University
2006 Sheridan Road
Evanston, Illinois 60208
USA
E-mail: yuanliao2008@u.northwestern.edu
         wjiang@northwestern.edu